\documentclass[leqno,12pt]{amsart}
\usepackage{amssymb,amscd,diagrams}

\calclayout

\newtheorem{theorem}{Theorem}
\newtheorem{corollary}{Corollary}
\newtheorem{proposition}{Proposition}
\newtheorem{lemma}{Lemma}

\newcommand{\cB}{{\mathcal B}}
\newcommand{\cC}{{\mathcal C}}
\newcommand{\cD}{{\mathcal D}}
\newcommand{\cF}{{\mathcal F}}
\newcommand{\cG}{{\mathcal G}}

\newcommand{\cL}{{\mathcal L}}
\newcommand{\cO}{{\mathcal O}}
\newcommand{\cX}{{\mathcal X}}

\newcommand{\mA}{{\mathbb A}}
\newcommand{\mC}{{\mathbb C}}
\newcommand{\mP}{{\mathbb P}}

\title{On orbit closures of Borel subgroups in spherical varieties}
\author{Michel~Brion}
\address{Universit\'e de Grenoble I\\
D\'epartement de Math\'ematiques\\
Institut Fourier, UMR 5582 du CNRS\\
38402 Saint-Martin d'H\`eres Cedex, France}
\email{Michel.Brion@ujf-grenoble.fr}

\date{}

\begin{document}

\begin{abstract}
Let $\cF$ be the flag variety of a complex semi-simple group $G$,
let $H$ be an algebraic subgroup of $G$ acting on $\cF$ with finitely
many orbits, and let $V$ be an $H$-orbit closure in $\cF$. Expanding the
cohomology class of $V$ in the basis of Schubert classes defines a
union $V_0$ of Schubert varieties in $\cF$ with positive
multiplicities. If $G$ is simply-laced, we show that these
multiplicites are equal to the same power of $2$. For arbitrary $G$, we 
show that $V_0$ is connected in codimension $1$. If moreover all
multiplicities are $1$, we show that the singularities of $V$ are
rational, and we construct a flat degeneration of $V$ to $V_0$. 
Thus, for any effective line bundle $L$ on $\cF$, the
restriction map $H^0(\cF,L)\to H^0(V,L)$ is surjective, and
$H^i(V,L)=0$ for $i\geq 1$.
\end{abstract}

\maketitle

\section*{Introduction}

Let $X$ be a spherical variety, that is, $X$ is a normal algebraic
variety endowed with an action of a connected reductive group $G$ such
that the set of orbits of a Borel subgroup $B$ in $X$ is finite. These 
$B$-orbits play an important role in the geometry and topology of 
$X$: they define a stratification by products of affine spaces with tori, 
and the Chow group of $X$ is generated by the classes of their closures. 
Moreover, the $B$-orbits in a spherical homogeneous space $G/H$, viewed 
as $H$-orbits in the flag variety $G/B$, are of importance in
representation theory.

\smallskip

The set $\cB(X)$ of $B$-orbit closures in $X$ is partially ordered by
inclusion. A weaker order $\preceq$ of $\cB(X)$ is defined by:
$Y\preceq Y'$ if there exists a sequence $(P_1,\ldots,P_n)$ of 
subgroups containing $B$ such that $Y'=P_1\cdots P_n Y$. In this
paper, we establish some properties of this weak order and its
associated graph, with applications to the geometry of $B$-orbit
closures.

\smallskip

Both orders are well known in the case where $X$ is the flag
variety of $G$. Then $\cB(X)$ identifies to the Weyl group $W$, and
the inclusion (resp. weak) order is the Bruhat-Chevalley (resp. left)
order, see e.g. \cite{H} 5.8. The $B$-orbit closures are the Schubert
varieties; their singularities are rational, in particular, they are
normal and Cohen-Macaulay.

\smallskip

Other important examples of homogeneous spherical varieties are
symmetric spaces. In this case, the inclusion and weak orders have been
studied in detail by Richardson and Springer \cite{RS1}, \cite{RS2},
\cite{Sp1}. But the geometry of $B$-orbit closures is far from being
fully understood; some of them are non-normal, see \cite{BE}. 

\smallskip

Returning to the general setting of spherical varieties, 
examples of $B$-orbit closures of arbitrary dimension and depth 1 are
given at the beginning of Section 3.  On the other hand, the
singularities of all $G$-orbit closures in a spherical $G$-variety are
rational, see e.g. \cite{BI}. A criterion for $B$-orbit closures to
have rational singularities will be formulated below, in terms of the
oriented graph $\Gamma(X)$ associated with the weak order.

\smallskip

For this, we endow $\Gamma(X)$ with additional data, as in \cite{RS1}:
each edge from $Y$ to $Y'$ is labeled by a simple root of $G$
corresponding to a minimal parabolic subgroup $P$ such that
$PY=Y'$. The degree of the associated morphism $P\times^B Y\to Y'$
being $1$ or $2$, this defines simple and double edges. There may be
several labeled edges with the same endpoints; but they are
simultaneously simple or double (Proposition \ref{type}).

\smallskip

For a spherical homogeneous space $G/H$, the cohomology classes of
$H$-orbit closures in $G/B$ can be read off the graph $\Gamma(G/H)$: 
each $H$-orbit closure $V$ in $G/B$ corresponds to a
$B$-orbit closure $Y$ in $X$. Consider an oriented path
$\gamma$ in $\Gamma(X)$, joining $Y$ to $X$. Denote by $D(\gamma)$ 
its number of double edges, and by $w(\gamma)$ the product in $W$ 
of the simple reflections associated with its labels. It turns out 
that $D(\gamma)$ depends only of $Y$ and $w(\gamma)$ 
(Lemma \ref{length}) and that we have in the cohomology ring of $G/B$:
$$
[V]=\sum_{w=w(\gamma)} 2^{D(\gamma)}\;[\overline{Bw_0wB}/B],
$$
the sum over the $w(\gamma)$ associated with all oriented paths from
$Y$ to $X$. Here $w_0$ denotes the longest element of $W$.

\smallskip

Thus, we are led to study oriented paths in $\Gamma(X)$ and their
associated Weyl group elements; this is the topic of Section 1. 
The main tool is a notion of neighbor paths that reduces several
questions to the case where $G$ has rank two. Using this, we show that
the union of Schubert varieties
$$
V_0=\bigcup_{w=w(\gamma)}\overline{Bw_0wB}/B
$$
is connected in codimension $1$ (Corollary \ref{flag}). If moreover
$G$ is simply-laced, then $D(\gamma)$ depends only on the
endpoints of $\gamma$ (Proposition \ref{constant}). As a consequence,
all coefficients of $[V]$ in the basis of Schubert classes are
equal. For symmetric spaces, the latter result is due to Richardson
and Springer \cite{Sp2}. It does not extend to multiply-laced groups,
see Example 3 in Section 1.

\smallskip

In Section 2, we analyze the intersections of $B$-orbit closures with
$G$-orbit closures in an important class of spherical varieties, the
(complete) regular $G$-varieties in the sense of Bifet, De Concini and
Procesi \cite{BDP}. This generalizes results of \cite{B1} \S 1 where
the intersections with closed $G$-orbits were described. Here the new 
ingredient is the construction of a ``slice'' $S_{Y,w}$ associated
with a $B$-orbit closure $Y$ in complete regular $X$, and with the
Weyl group element $w$ defined by an oriented path from $Y$ to
$X$. The $S_{Y,w}$ are toric varieties; each oriented path $\gamma$ in
$\Gamma(X)$ defines a finite surjective morphism between ``slices'' of
its endpoints, of degree $2^{D(\gamma)}$. If the target of $\gamma$ is
$X$, then the intersection multiplicities of $Y$ with all $G$-orbit
closures that meet $S_{Y,w}$ turn out to be divisors of $2^{D(\gamma)}$.
Moreover, given a $G$-orbit closure $X'$ and an irreducible component
$Y'$ of $Y\cap X'$, there exists a ``slice'' meeting $Y'$ (Theorem
\ref{indices}.)

\smallskip

This distinguishes the $B$-orbit closures $Y$ such that all
oriented paths in $\Gamma(X)$ with source $Y$ contain simple edges
only; we call them multiplicity-free. In a regular variety, any 
irreducible component of the intersection of multiplicity-free $Y$
with a $G$-orbit closure is multiplicity-free as well, and the
corresponding intersection multiplicity equals $1$ (Corollary
\ref{components}.)

\smallskip

Section 3 contains our main result: the singularities of any
multiplicity-free $B$-orbit closure $Y$ in a spherical variety $X$ are
rational, if $X$ contains no fixed points of simple normal subgroups
of $G$ of type $G_2$, $F_4$ and $E_8$ (Theorem \ref{normal}; its
technical assumption is used in one of the reduction steps of the
proof, but the statement should hold in full generality.) The proof
goes by decreasing induction on $Y$, like Seshadri's proof of
normality of Schubert varieties \cite{Se}. This result applies, e.g.,
to regular $G$-varieties; for them, we show that the
scheme-theoretical intersection of $Y$ with any $G$-orbit closure is
reduced.

\smallskip

For a $H$-orbit closure $V$ in $G/B$, the corresponding $B$-orbit
closure $Y$ is multiplicity-free if and only if $[V]=[V_0]$. In that
case, we construct a flat degeneration of $V$ to $V_0$, where the latter
is viewed as a reduced subscheme of $G/B$ (Corollary \ref{flat}). 
Thus, the equality $[V]=[V_0]$ holds in the Grothendieck group of
$G/B$ as well. As another consequence, the restriction map
$H^0(G/B,L)\to H^0(V,L)$ is surjective for any effective line bundle
$L$ on $G/B$; moreover, the higher cohomology groups $H^i(V,L)$ vanish
for $i\geq 1$ (Corollary \ref{surj}.)
Applied to symmetric spaces and combined with Theorem B of \cite{BE},
the latter result implies a version of the 
Parthasaraty-Ranga Rao-Varadarajan conjecture, see \cite{BE} \S 6. 
It extends to certain smooth $H$-orbit closures, but not to all of
them, see the example in \cite{BH} 4.3. In fact, surjectivity of all
restriction maps for spherical $G/H$ is equivalent to
multiplicity-freeness of all $H$-orbit closures in $G/B$ 
(Proposition \ref{nonsurj}.) 

\smallskip

In Section 4, we relate our approach to work of Knop
\cite{K1}, \cite{K2}. He defined an action of $W$ on $\cB(X)$ such
that the $W$-conjugates of the maximal element $X$ are the orbit
closures of maximal rank (in the sense of \cite{K2}). Moreover, the
isotropy group $W_{(X)}$ is closely related to the ``Weyl group of
$X$'', as defined in \cite{K1}. It is easy to see that all orbit
closures of maximal rank are multiplicity-free, and hence their
singularities are rational if $X$ is regular. In that case, we
describe the intersections of $B$-orbit closures of maximal rank with
$G$-orbit closures, in terms of $W$ and $W_{(X)}$ (Proposition
\ref{minimal}.)

\smallskip

This implies two results on the position of $W_{(X)}$ in $W$: firstly,
all elements of $W$ of minimal length in a given $W_{(X)}$-coset have
the same length. Secondly, $W_{(X)}$ is generated by reflections or
products of two commuting reflections of $W$. This gives a
simple proof of the fact that the Weyl group of $X$ is generated by
reflections \cite{K1}.

\smallskip

A remarkable example of a spherical homogeneous space where all
orbit closures of a Borel subgroup have maximal rank is the group $G$
viewed as a homogeneous space under $G\times G$. If moreover $G$ is
adjoint, then it has a canonical $G\times G$-equivariant completion
${\bf X}$. It is proved in \cite{BP} that the $B\times B$-orbit 
closures in ${\bf X}$ are normal, and that their intersections
are reduced. This follows from the fact that ${\bf X}$ is Frobenius 
split compatibly with all $B\times B$-orbit closures.

\smallskip

It is tempting to generalize this to any spherical variety $X$. By
\cite{BI}, $X$ is Frobenius split compatibly with all $G$-orbit
closures. But this does not extend to $B$-orbit closures, since their
intersections may be not reduced. This happens, e.g., for the space 
of all symmetric $n\times n$ matrices of rank $n$, that is, the
symmetric space ${\rm GL}(n)/{\rm O}(n)$: consider the subvarieties
$(a_{11}=0)$ and $(a_{11}a_{22}-a_{12}^2=0)$. On the other hand, many
$B$-orbit closures in that space are not normal for $n\geq 5$, see
\cite{P}.

\smallskip

So the present paper generalizes part of the results of \cite{BP} 
to all spherical varieties, by other methods. It raises many further
questions, e.g., is it true that the normalization of any $B$-orbit
closure in a spherical variety has rational singularities ? And do
our results extend to positive characteristics (the proof of Theorem
\ref{normal} uses an equivariant resolution of singularities) ?

\medskip

\noindent
{\sl Acknowledgements}. I thank Peter Littelmann, Laurent Manivel,
Olivier Mathieu, St\'ephane Pin, Patrick Polo and Tonny Springer for
useful discussions or e-mail exchanges.
\medskip

\noindent
{\sl Notation.} Let $G$ be a complex connected reductive algebraic
group. Let $B$ be a Borel subgroup of $G$, with unipotent radical
$U$. Let $T$ be a maximal torus of $B$, with Weyl group $W$. Let $\cX$
be the character group of $B$; we identify $\cX$ with the character
group of $T$, and we choose a $W$-invariant scalar product on
$\cX$. Let $\Phi$ be the root system of $(G,T)$, with the subset
$\Phi^+$ of positive roots defined by $B$, and its subset $\Delta$ of
simple roots. 

For $\alpha\in\Delta$, let $s_{\alpha}\in W$ be the corresponding
simple reflection, and let $P_{\alpha}=B\cup Bs_{\alpha}B$ be the
corresponding minimal parabolic subgroup. For any subset $I$ of
$\Delta$, let $P_I$ be the subgroup of $G$ generated by the
$P_{\alpha}$, $\alpha\in I$. The map $I\mapsto P_I$ is a bijection
from subsets of $\Delta$ to subgroups of $G$ containing $B$, that is,
to standard parabolic subgroups of $G$.

Let $L_I$ be the Levi subgroup of $P_I$ that contains $T$; let
$\Phi_I$ be the root system of $(L_I,T)$, with Weyl group $W_I$. We
denote by $\ell$ the length function on $W$ and by $W^I$ the set of
all $w\in W$ such that $\ell(ws_{\alpha})=\ell(w)+1$ for all
$\alpha\in I$ (this amounts to: $w(I)\subseteq \Phi^+$). Then $W^I$
is a system of representatives of the set of right cosets $W/W_I$. 

\section{The weak order and its graph}

In the sequel, we denote by $X$ a complex spherical $G$-variety and 
by $\cB(X)$ the set of $B$-orbit closures in $X$.
One associates to a given $Y\in\cB(X)$ several combinatorial
invariants, see \cite{K2}: The {\sl character group} $\cX(Y)$
is the set of all characters of $B$ that arise as weights of
eigenvectors of $B$ in the function field $\mC(Y)$. Then $\cX(Y)$ is a
free abelian group of finite rank $r(Y)$, the {\sl rank} of $Y$.  

Let $Y^0$ be the open $B$-orbit in $Y$ and let $P(Y)$ be the  
set of all $g\in G$ such that $gY^0=Y^0$; then $P(Y)$ is a standard
parabolic subgroup of $G$. Let $L(Y)$ be its Levi subgroup that
contains $T$ and let $\Delta(Y)$ be the corresponding subset of
$\Delta$: the set of {\sl simple roots of} $Y$.

We note some easy properties of these invariants.

\begin{lemma}\label{prelim}
(i) $\cX(Y)$ is isomorphic to the quotient of the group of invertible
regular functions on $Y^0$, by the subgroup of constant non-zero
functions.

\noindent
(ii) The derived subgroup $[L(Y),L(Y)]$ fixes a point of $Y^0$.

\noindent
(iii) The group $W_{\Delta(Y)}$ fixes pointwise
$\cX(Y)$. Equivalently, any simple root of $Y$ is orthogonal to
$\cX(Y)$.
\end{lemma}

\begin{proof} (i) Let $f$ be an eigenvector of $B$ in $\mC(Y)$ with
weight $\chi(f)$. Then $f$ restricts to an invertible regular function
on $Y^0$, and is uniquely determined by $\chi(f)$ up to a
constant. Conversely, let $f$ be an invertible regular function on
the $B$-orbit $Y^0$. Then $f$ pulls back to an invertible regular
function on $B$, that is, to a scalar multiple of a
character of $B$. Thus, $f$ is an eigenvector of $B$ in $\mC(Y)$.

\noindent
(ii) Choose $y\in Y^0$. Let $B_y$ (resp. $P(Y)_y$) be
the isotropy group of $y$ in $B$ (resp. $P(Y)$). Since $Y^0=By=P(Y)y$,
we have $P(Y)=BP(Y)_y$. Thus, $P(Y)_y$ acts transitively on $P(Y)/B$,
the flag variety of $P(Y)$. Using e.g. \cite{Dem}, it follows that
$P(Y)_y$ contains a maximal connected semisimple subgroup of $P(Y)$,
that is, a conjugate of $[L(Y),L(Y)]$.

\noindent
(iii) follows from \cite{K2} Lemma 3.2; it can be deduced from (ii) as
well.
\end{proof}

Let $\cD(X)$ be the subset of $\cB(X)$ consisting of irreducible
$B$-stable divisors that are not $G$-stable. The elements of $\cD(X)$
are called {\sl colors}; they play an important role in the
classification of spherical embeddings, see \cite{K0}. They also allow
to describe the parabolic subgroups associated with $G$-orbit
closures:

\begin{lemma}\label{parab}
Let $Y$ be the closure of a $G$-orbit in $X$ and let $\cD_Y(X)$ be the
set of all colors that contain $Y$. Then $P(Y)$ is the set of all
$g\in G$ such that $gD=D$ for any $D\in\cD(X)-\cD_Y(X)$. Moreover, 
there exists $y\in Y^0$ fixed by $[L(Y),L(Y)]$, such that the map 
$R_u(P(Y))\times Ty\to Y^0, (g,x)\mapsto gx$ is an isomorphism. Then
the dimension of $Ty$ equals the rank of $Y$.
\end{lemma}

\begin{proof}
Let $X_0$ be the complement in $X$ of the union of all irreducible
$B$-stable divisors that do not contain $Y$. Then $X_0$ is an open
affine $B$-stable subset of $X$, and $X_0\cap Y$ equals $Y^0$; see
\cite{K0} Theorem 3.1. Let $Q$ be the stabilizer of $X_0$ in $G$, then
$Q$ consists of all $g\in G$ such that $gD=D$ for all 
$D\in\cD(X)-\cD_Y(X)$. Clearly, $Q$ is a standard parabolic subgroup,
contained in $P(Y)$. It follows that $R_u(P(Y))\subseteq R_u(Q)$.

Let $M$ be the standard Levi subgroup of $Q$. By \cite{K1} 2.3 and
2.4, there exists a closed $M$-stable subvariety $S$ of $X_0$ such
that the product map $R_u(Q)\times S\to X_0$ is an isomorphism;
moreover, $[M,M]$ acts trivially on $S\cap Y^0$. In particular, for
any $y\in S\cap Y^0$, the product map $R_u(Q)\times Ty\to Y^0$ is an 
isomorphism. Since $R_u(Q)=R_u(P(Y)) (R_u(Q)\cap[L(Y),L(Y)])$ and
since $[L(Y),L(Y)]$ fixes points of $Y^0$, it follows that
$R_u(Q)=R_u(P(Y))$, whence $Q=P(Y)$. Moreover, the character group of
$Y$ is isomorphic to that of the torus $Ty\cong T/T_Y$, whence
$r(Y)=\dim(Ty)$.
\end{proof}

This description of $Y^0$ as a product of a unipotent group with a
torus will be generalized in Section 4 to all $B$-orbits of maximal
rank.

Returning to arbitrary $B$-orbit closures, let $Y,Y'\in\cB(X)$ and let
$\alpha\in\Delta$. We say that $\alpha$ {\sl raises $Y$ to} $Y'$ if
$Y'=P_{\alpha}Y\neq Y$. Let then
$$
f_{Y,\alpha}:P_{\alpha}\times^B Y\to P_{\alpha}/B
$$ 
be the homogeneous bundle with fiber the $B$-variety $Y$ and basis
$P_{\alpha}/B$ (isomorphic to projective line.) The map
$P_{\alpha}\times Y\to X$,$(p,y)\mapsto py$ 
factors through a proper morphism
$$
\pi_{Y,\alpha}:P_{\alpha}\times^B Y\to 
Y'=P_{\alpha}Y
$$
that restricts to a finite morphism
$P_{\alpha}\times^B Y^0\to P_{\alpha}Y^0$.

By \cite{RS1} or \cite{K2} Lemma 3.2, one of the following three cases
occurs.

\noindent
$\bullet$ Type $U$: $P_{\alpha}Y^0=Y^{'0}\cup Y^0$ and $\pi_{Y,\alpha}$
is birational. Then $\cX(Y')=s_{\alpha}\cX(Y)$; thus, $r(Y')=r(Y)$.

\noindent
$\bullet$ Type $T$: $P_{\alpha}Y^0=Y^{'0}\cup Y^0\cup Y^0_-$ for some 
$Y_-\in\cB(X)$ of the same dimension as $Y$, and $\pi_{Y,\alpha}$ is
birational. Then $r(Y)=r(Y_-)=r(Y')-1$. 

\noindent
$\bullet$ Type $N$: $P_{\alpha}Y^0=Y^{'0}\cup Y^0$ and $\pi_{Y,\alpha}$
has degree 2. Then $r(Y)=r(Y')-1$.

In particular, $r(Y)\leq r(P_{\alpha}Y)$ with equality if and only if
$\alpha$ has type $U$ .

Our notation for types differs from that in \cite{RS1} and \cite{K2};
it can be explained as follows. Choose $y\in Y^0$ with isotropy group
$(P_{\alpha})_y$ in $P_{\alpha}$. Then $(P_{\alpha})_y$ acts on
$P_{\alpha}/B\cong\mP^1$ with finitely many orbits, for $B$ acts on
$P_{\alpha}Y^0\cong P_{\alpha}/(P_{\alpha})_y$ with finitely many
orbits. By \cite{RS1} or \cite{K2}, the image of $(P_{\alpha})_y$ in 
${\rm Aut}(P_{\alpha}/B)\cong{\rm PGL}(2)$ is a torus (resp. the
normalizer of a torus) in type $T$ (resp. $N$); in type $U$, this image 
contains a non-trivial unipotent normal subgroup.
\medskip

\noindent
{\sl Definition.} Let $\Gamma(X)$ be the oriented graph with vertices
the elements of $\cB(X)$ and edges labeled by $\Delta$, where $Y$ is
joined to $Y'$ by an edge of label $\alpha$ if that simple root raises
$Y$ to $Y'$. This edge is simple (resp. double) if $\pi_{Y,\alpha}$ has
degree 1 (resp. 2.) The partial order $\preceq$ on $\cB(X)$ with oriented
graph $\Gamma(X)$ will be called the {\sl weak order}.
\medskip

Observe that the dimension and rank functions are
compatible with $\preceq$. 
We shall see that $Y,Y'\in\cB(X)$ satisfy $Y\preceq Y'$ if and only if
there exists $w\in W$ such that $Y'$ equals the closure
$\overline{BwY}$ (Corollary \ref{weak}.)

In the case where $X=G/P$ where $P$ is a parabolic subgroup of $G$,
the rank function is zero. Thus, all edges are of type $U$; in
particular, they are simple.

Here is another example, where double edges occur.
\medskip

\noindent
{\sl Example 1.} Let $G={\rm GL}(3)$ with simple roots $\alpha$ and
$\beta$. Let $H$ be the subgroup of $G$ consisting of matrices of the
form
$$
\left(\begin{matrix}
*&0&*\cr 0&*&*\cr 0&0&*\cr
\end{matrix}\right)~{\rm or}
\left(\begin{matrix}
0&*&*\cr *&0&*\cr 0&0&*\cr
\end{matrix}\right)
$$
and let $X=G/H$. It is easy to see that $X$ is spherical of rank
one and that $\Gamma(X)$ is as follows:

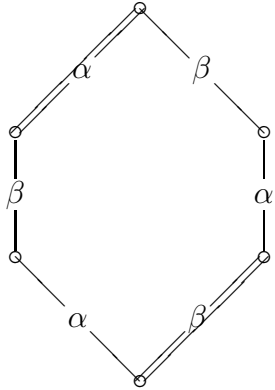
\begin{diagram}[abut,height=2em,width=2em]
&&\circ&&&&\\
&\ldLine\;\ldLine~\alpha&&\rdLine~\beta\\
\circ&&&&\circ&\\
\dLine~\beta&&&&\dLine~\alpha\\
\circ&&&&\circ\\
&\rdLine~\alpha&&\ldLine~\beta\;\ldLine\\
&&\circ&&\\
\end{diagram}

Observe that $\Gamma(X)$ is the same as $\Gamma(G/B)$, except for
double edges. But the geometry of $B$-orbit closures is very different
in both cases: all of them are smooth in $G/B$ (the flag variety of
$\mP^2$), whereas $X$ contains a $B$-stable divisor that is singular
in codimension $1$. 

Specifically, let $Z$ be the closed $B$-orbit in $G/H$. We claim that
$Y=P_{\beta}P_{\alpha}Z$ is singular along $P_{\beta}Z$. Indeed, the
morphism $\pi:P_{\beta}\times^B P_{\alpha}Z\to Y$
is birational, and $\pi^{-1}(P_{\beta}Z)$ equals
$P_{\beta}\times^B Z$. But the restriction 
$P_{\beta}\times^B Z\to P_{\beta}Z$ has degree two. Now our claim
follows from Zariski's main theorem.

One checks that $r(P_{\beta}Z)=1$, whereas $r(Y)=0$. Thus, the
rank function is not compatible with the inclusion order.
\medskip

Returning to the general situation, observe that $GY$ is the closure
of a $G$-orbit for any $Y\in\cB(X)$. Moreover, $Y$ is 
the source of an oriented path in $\Gamma(X)$ with target $GY$, since 
the group $G$ is generated by the $P_{\alpha}$, $\alpha\in\Delta$. 
By \cite{K2} Corollary 2.4, we have $r(GY)\leq r(X)$, so that 
$r(Y)\leq r(X)$. It also follows that each connected component of 
$\Gamma(X)$ contains a unique $G$-orbit closure. 

The simple roots of $Y$ are determined by $\Gamma(X)$:
indeed, $\alpha\in\Delta$ is not in $\Delta(Y)$ if and only if
$\alpha$ is the label of an edge with endpoint $Y$. 
Similarly, if $\alpha$ raises $Y$ then its type is determined by
$\Gamma(X)$: it is $U$ (resp. $N$) if there is a unique edge of label
$\alpha$ and target $P_{\alpha}Y$ and this edge is simple
(resp. double); and it is $T$ if there are two such edges.
It follows that the ranks of $B$-orbit closures are 
determined by $\Gamma(X)$ and the ranks of $G$-orbit closures.

There is no restriction on the number of edges in $\Gamma(X)$ with
prescribed endpoints, as shown by the example below suggested by
D.~Luna. But we shall see that all such edges have the same type.
\medskip

\noindent
{\sl Example 2.} Let $n$ be a positive integer. Let 
$G={\rm SL}(2)\times\cdots\times{\rm SL}(2)$ ($n$ terms) and let $H$
be the subgroup of $G$ consisting of those $n$-tuples
$$
\left(\begin{matrix}t&u_1\cr 0&t^{-1}\cr\end{matrix}\right),
\ldots,
\left(\begin{matrix}t&u_n\cr 0&t^{-1}\cr\end{matrix}\right)
$$
where $t\in\mC^*$, $u_1,\ldots,u_n\in\mC$ and $u_1+\cdots+u_n=0$.
One checks that $G/H$ is spherical; the open $H$-orbit in 
$G/B\cong {\mP}^1\times\cdots\times{\mP}^1$ ($n$ terms) consists of 
those $(z_1,\ldots,z_n)$ such that $z_i\neq\infty$
for all $i$, and that $z_1+\cdots+z_n\neq 0$. Let $Y$ be the $B$-stable 
hypersurface in $G/H$ corresponding to the $H$-stable hypersurface 
$(z_1+\cdots+z_n=0)$ in $G/B$. One checks that $Y$ is irreducible and
raised to $G/H$ by all simple roots of $G$ (there are $n$ of
them). Thus, $Y$ is joined to $G/H$ by $n$ edges of type $U$.
\medskip

\begin{proposition}\label{type}
Let $Y,Y'\in\cB(X)$ and let $\alpha$, $\beta$ be distinct
simple roots raising $Y$ to $Y'$. Then either $\alpha$, $\beta$
are orthogonal and both of type $U$, or they are both of type $T$.
\end{proposition}

\begin{proof} We begin with two lemmas that reduce the ``local'' study
of $\Gamma(X)$ to simpler situations.

Let $Y\in\cB(X)$ and let $P=P_I$ be a standard parabolic subgroup of
$G$, with radical $R(P)$. Let $\cB(P,Y)$ be the set of all closures in
$X$ of $B$-orbits in $PY^0$; in other words, $\cB(P,Y)$ is the set of
all $Z\in\cB(X)$ such that $PZ=PY$. Let $\Gamma(P,Y)$ be the oriented
graph with set of vertices $\cB(P,Y)$, and with edges those edges of
$\Gamma(X)$ that have both endpoints in $\cB(P,Y)$ and labels in $I$.

\begin{lemma}\label{reduc}
The quotient $PY^0/R(P)$ is a $P/R(P)$-homogeneous spherical variety
with graph $\Gamma(P,Y)$.
\end{lemma}

\begin{proof} Since $PY^0$ is a unique $P$-orbit and $R(P)$ is a normal
subgroup of $P$ contained in $B$, the quotient $PY^0/R(P)$ exists and
is homogeneous under $P/R(P)$; moreover, any $B/R(P)$-orbit in
$PY^0/R(P)$ pulls back to a unique $B$-orbit in $PY^0$. Let $\cO$ be a
$B$-orbit in $PY^0$ and let $\alpha\in I$. Then $R(P_{\alpha})$
contains $R(P)$, the square
$$
\begin{matrix}
P_{\alpha}\times^B\cO&\to&P_{\alpha}\cO\cr
\downarrow&&\downarrow\cr
P_{\alpha}\times^B\cO/R(P)&\to&P_{\alpha}\cO/R(P)\cr
\end{matrix}
$$
is cartesian, and the map 
$P_{\alpha}\times^B\cO/R(P)\to P_{\alpha}/R(P)\times^{B/R(P)}\cO/R(P)$
is an isomorphism. Thus, the type is preserved under pull back.
\end{proof}

Assume now that $X$ is homogeneous under $G$; write then $X=G/H$.
Let $H'$ be a closed subgroup of the normalizer $N_G(H)$ such
that $H'$ contains $H$, and that the quotient $H'/H$ is connected. Let
$Z(G)$ be the center of $G$. Let $X'=G/H'Z(G)$, a homogeneous
spherical variety under the adjoint group $G/Z(G)$. The natural
$G$-equivariant map $p:X\to X'$ is the quotient by the right action of
$H'Z(G)$ on $G/H$. 

\begin{lemma}\label{nor}
The pull-back under $p$ of any $B$-orbit in $X'$ is a unique $B$-orbit 
in $X$. This defines an isomorphism of $\Gamma(X')$ onto $\Gamma(X)$.
\end{lemma}

\begin{proof}
The first assertion follows from \cite{B2} Proposition 2.2 (iii). The
second assertion is checked as in the proof of Lemma \ref{reduc}.
\end{proof}

\begin{lemma}\label{orth}
Let $Y\in\cB(X)$, $Y\neq X$, and let $\alpha\in\Delta$. If
$P_{\alpha}Y^0=X$ then $\alpha$ is orthogonal to $\Delta-\{\alpha\}$,
and the derived subgroup of $L_{\Delta -\{\alpha\}}$ fixes
pointwise $X$.
\end{lemma}

\begin{proof} Let $H$ be the isotropy group in $G$ of a point of
$Y^0$. Since $P_{\alpha}Y^0=X$, we have $P_{\alpha}H=G$. Equivalently, 
the map $H/P_{\alpha}\cap H\to G/P_{\alpha}$ is an isomorphism. But since 
$Y\neq X$, we have $Y^0\neq P_{\alpha}Y^0$, so that the image of
$P_{\alpha}\cap H$ in $P_{\alpha}/R(P_{\alpha})\cong{\rm PGL}(2)$ is a 
proper subgroup. It follows that $(P_{\alpha}\cap H)^0$ is solvable.
Thus, $H/P_{\alpha}\cap H$ is the flag variety of $H^0$. Now the 
connected automorphism group of this flag variety is the quotient of
$H^0/R(H^0)$ by its center. On the other hand, the connected
automorphism group of $G/P_{\alpha}$ is $G/Z(G)$ if $\alpha$ is not
orthogonal to $\Delta-\{\alpha\}$ (this follows e.g. from
\cite{Dem}.) In this case, we have $G=Z(G)H^0$ so that
$G/H$ is a unique $B$-orbit, a contradiction. Thus, $G/Z(G)$ is
the product of $L_{\alpha}/Z(L_{\alpha})$ with 
$L_{\Delta -\{\alpha\}}/Z(L_{\Delta -\{\alpha\}})$, and the map 
$L_{\Delta-\{\alpha\}}/B\cap L_{\Delta-\{\alpha\}}\to G/P_{\alpha}$
is an isomorphism. It follows that the derived subgroup of 
$L_{\Delta -\{\alpha\}}$ is contained in $H$.
\end{proof}

We now prove Proposition \ref{type}. Applying Lemma \ref{reduc} to
$Y'$ and $P_{\alpha,\beta}$, we may assume that $Y'=X=G/H$ for some 
subgroup $H$ of $G$ and that $\Delta=\{\alpha,\beta\}$. 

If $\alpha$ has type $U$, then $r(Y)=r(X)$ whence $\beta$ has type $U$
as well. We claim that $\cB(X)$ consists of $Y$ and $X$. Indeed, if
$Z\in\cB(X)$ and $Z\neq X$, then $Z$ is connected to $X$ by an
oriented path in $\Gamma(X)$. Let $Z'$ be the source of the top edge
of this path. That edge cannot have $Y$ as its target (otherwise $Y$
would be stable under $P_{\alpha}$ or $P_{\beta}$); thus, it raises
$Z'$ to $X$. Since $\alpha$ and $\beta$ have type $U$, it follows that
$Z'=Y$, whence $Z=Y$. Thus, $P_{\alpha}Y^0=X$; then $\alpha$ and
$\beta$ are orthogonal by Lemma \ref{orth}. 

If $\alpha$ has type $N$, then $r(Y)=r(X)-1$, whence $\beta$ has type
$N$ or $T$. In the former case, we see as above that
$X=P_{\alpha}Y^0=P_{\beta}Y^0$. Thus, $\alpha$ and $\beta$ are
orthogonal by Lemma \ref{orth}. Using Lemma \ref{nor}, we may assume
that $G={\rm PGL}(2)\times{\rm PGL}(2)$ and that $H$ contains a copy
of ${\rm PGL}(2)$. Then $H$ is conjugate to ${\rm PGL}(2)$ embedded 
diagonally in $G$. But then both $\alpha$ and $\beta$ have type $T$, 
a contradiction.

If $\alpha$ has type $N$ and $\beta$ has type $T$, then there exists
$y\in Y^0$ such that $(P_{\beta})_y$ is contained in $R(P_{\beta})T$.
Since the homogeneous spaces $P_{\beta}/R(P_{\beta})T$ and
$R(P_{\beta})T/(P_{\beta})_y$ are affine, the same holds for
$P_{\beta}/(P_{\beta})_y\cong P_{\beta}Y^0$. It follows that 
$X - P_{\beta}Y^0$ is pure of codimension $1$ in $X$. But
$P_{\beta}Y^0$ meets both $B$-orbits of codimension $1$ in $X$, so
that $P_{\beta}Y^0=X$. This case is excluded as above. Thus, type $N$
does not occur.
\end{proof}

We next study oriented paths in $\Gamma(X)$.
Let $\gamma$ be such a path, with source $Y$ and target $Y'$. Let
$(\alpha_1,\alpha_2,\ldots,\alpha_{\ell})$ be the sequence of labels
of edges of $\gamma$, where $\ell=\ell(\gamma)$ is the length of the
path. Let $\ell_U(\gamma)$ (resp. $\ell_T(\gamma)$, $\ell_N(\gamma)$)
be the number of edges of type $U$ (resp. $T$, $N$) in $\gamma$. Then
$$
\ell_U(\gamma)+\ell_T(\gamma)+\ell_N(\gamma)=\ell(\gamma)=
\dim(Y')-\dim(Y).
$$
Define an element $w(\gamma)$ of $W$ by
$w(\gamma)=s_{\alpha_{\ell}}\cdots s_{\alpha_2}s_{\alpha_1}$.

\begin{lemma}\label{length}
(i) $(s_{\alpha_\ell},\ldots,s_{\alpha_2},s_{\alpha_1})$ is a
reduced decomposition of $w(\gamma)$; equivalently, 
$\ell(w(\gamma))=\ell$. 

\noindent
(ii) $\ell_T(\gamma)+\ell_N(\gamma)=r(Y')-r(Y)$. In
particular, $\ell_T(\gamma)+\ell_N(\gamma)$ and $\ell_U(\gamma)$
depend only on the endpoints of $\gamma$.

\noindent
(iii) The morphism $G\times^B Y\to X:(g,y)B\to gy$
restricts to a morphism 
$\overline{Bw(\gamma)B}\times^B Y\to Y'$
that is surjective and generically finite
of degree $2^{\ell_N(\gamma)}$. In particular, $\ell_T(\gamma)$ and
$\ell_N(\gamma)$ depend only on the endpoints of $\gamma$ and on
$w(\gamma)$. Moreover, $w(\gamma)$ is in $W^{\Delta(Y)}$, and
$w(\gamma)^{-1}$ is in $W^{\Delta(Y')}$.

\noindent
(iv) If the stabilizer in $G$ of a point of $Y^0$ is contained in a
Borel subgroup of $G$ (e.g., if $X=G/H$ where $H$ is connected and
solvable), then $\ell_N(\gamma)=0$ so that $\ell_T(\gamma)$
depends only on the endpoints of $\gamma$.
\end{lemma}

\begin{proof} (i) Observe that $Bs_{\alpha_1}Y$ is dense in
$P_{\alpha_1}Y$, as $P_{\alpha_1}$ raises $Y$. By induction, it
follows that 
$Bs_{\alpha_{\ell}}B\cdots s_{\alpha_2}Bs_{\alpha_1}Y$ 
is dense in $Y'$. Because $\dim(Y')=\dim(Y)+\ell$, we must have 
$\dim(\overline{Bs_{\alpha_{\ell}}B\cdots 
s_{\alpha_2}Bs_{\alpha_1}B}/B)=\ell$,
whence $\ell(s_{\alpha_{\ell}}\cdots s_{\alpha_2}s_{\alpha_1})=\ell$.

(ii) follows from the fact that $r(Y')=r(Y)$ (resp. $r(Y)+1$) if $Y$
is the source of an edge with target $Y'$ and type $U$ (resp. $T$,
$N$).

(iii) By (i), the product maps
$$
P_{\alpha_i}\times^B\cdots \times^B P_{\alpha_2}\times^B P_{\alpha_1}
\to\overline{Bs_{\alpha_i}\cdots s_{\alpha_2}s_{\alpha_1}B}
$$
are birational for $1\leq i\leq \ell$. It follows that the morphism 
$\overline{Bw(\gamma)B}\times^B Y\to X$ has image $Y'$; moreover, its
degree is the product of the degrees of the 
$$
\pi_i:P_{\alpha_i}\times^B (P_{\alpha_{i-1}}\cdots P_{\alpha_1}Y)
\to P_{\alpha_i}P_{\alpha_{i-1}}\cdots P_{\alpha_1}Y,
$$
that is, $2^{\ell_N(\gamma)}$.

Let $w=w(\gamma)$. We show that $w^{-1}\in W^{\Delta(Y')}$.
Otherwise, there exists $\alpha\in\Delta(Y')$ such that 
$\ell(s_{\alpha}w)=\ell(w)-1$. Thus, 
$BwB=Bs_{\alpha}Bs_{\alpha}wB$, and 
$Y'=\overline{BwY}=\overline{Bs_{\alpha}Bs_{\alpha}wY}$. Let
$Y''=\overline{Bs_{\alpha}wY}$, then $\alpha$ raises $Y''$ to
$Y'$. This contradicts the assumption that $\alpha\in\Delta(Y')$.
A similar argument shows that $w\in W^{\Delta(Y)}$.

(iv) If $\ell_N(\gamma)>0$, then there exists a point $x\in GY^0$,
a simple root $\alpha$ and a surjective group homomorphism 
$(P_{\alpha})_x\to N$ where $N$ is the normalizer of a torus in 
${\rm PGL}(2)$. Since $N$ consists of semisimple elements, it is a
quotient of $(P_{\alpha})_x/R_u(P_{\alpha})_x$. By assumption, the
latter is isomorphic to a subgroup of $B/U=T$. Thus, $N$ is abelian, 
a contradiction.

\end{proof}

\begin{corollary}\label{weak}
Let $Y,Y'\in\cB(X)$, then $Y\preceq Y'$ if and only if there exists 
$w\in W$ such that $Y'=\overline{BwY}$.
\end{corollary}

\begin{proof}
Recall that $\overline{BwB}$ (closure in $G$) is a product of minimal
parabolic subgroups. Thus,
$Y\preceq\overline{BwB}Y=\overline{BwY}$. The converse has just been
proved.
\end{proof}

For later use, we study the behavior of $\Gamma(X)$ under parabolic
induction in the following sense (see \cite{B1} 1.2.) Let $P=P_I$ be a
standard parabolic subgroup with Levi subgroup $L=L_I$ and let $X'$ be
a spherical $L$-variety, then the induced variety is $X=G\times^P X'$ 
where $P$ acts on $X'$ through its quotient $P/R_u(P)$, isomorphic to
$L$. In other words, $X$ is the total space of the homogeneous bundle
over $G/P$ with fiber $X'$. By [{\sl loc. cit.}], each $Y\in\cB(X)$
can be written uniquely as $\overline{BwY'}$ for $w\in W^I$ and
$Y'\in\cB(X')$; then $r(Y)=r(Y')$. We thus identify $\cB(X)$ to
$W^I\times \cB(X')$. The next result describes the edges of
$\Gamma(X)$ in terms of those of $\Gamma(X')$.

\begin{lemma}\label{induced}
Let $\alpha\in\Delta$, $w\in W^I$ and $Y'\in\cB(X')$; let 
$\beta=w^{-1}(\alpha)$. Then the edges of $\Gamma(X)$ with source
$(w,Y')$ and label $\alpha$ are as follows: 

\noindent
(i) If $\beta\in\Phi^+ - I$, join $(w,Y')$ to
$(s_{\alpha}w,Y')$ by an edge of type $U$.

\noindent
(ii) If $\beta\in I$ and $P_{\beta}\cap L$ raises $Y'$, join $(w,Y')$
to $(w,(P_{\beta}\cap L)Y')$ by an edge of the same type as the edge
from $Y'$ to $(P_{\beta}\cap L)Y'$.
\end{lemma}

\begin{proof}
Since $w\in W^I$, we have $s_{\alpha}w\in W^I$ if and only if
$\beta\notin I$. In that case, $P_{\alpha}$ raises $Y$
if and only if $\ell(s_{\alpha}w)=\ell(w)+1$, that is,
$\beta\in\Phi^+$. Then $P_{\alpha}Y=\overline{Bs_{\alpha}wY'}$
and the map $\pi_{Y,\alpha}$ is the pull-back of
$\pi_{\overline{BwP}/P,\alpha}$ under the map 
$\overline{BwY'}\to \overline{BwP}/P$. This yields case (i).

But if $\beta\in I$, then $s_{\alpha}w=ws_{\beta}$ has length
$\ell(w)+1$, so that
$$
P_{\alpha}Y=\overline{Bs_{\alpha}BwY'}=\overline{Bs_{\alpha}wY'}
=\overline{Bws_{\beta}Y'}=\overline{BwBs_{\beta}Y'}
=\overline{Bw(P_{\beta}\cap L)Y'}.
$$
Thus, $P_{\alpha}$ raises $Y$ if and only if $P_{\beta}\cap L$ raises
$Y'$. Then, as $s_{\alpha}w=ws_{\beta}$, we can join $Y'$ to
$P_{\alpha}Y$ by two paths: one beginning with $\ell(w)$ edges of type
$U$ followed by an edge from $Y$ to $P_{\alpha}Y$, and another one
beginning with an edge from $Y'$ to $(P_{\beta}\cap L)Y'$ followed by
$\ell(w)$ edges of type $U$. Using Lemma \ref{length}, this yields
case (ii). 
\end{proof}

For instance, Example 1 is obtained from ${\rm SL}(2)/N$ by parabolic
induction. 

Returning to the case where $X$ is an arbitrary spherical $G$-variety,
we shall see that the numbers $\ell_T(\gamma)$ and $\ell_N(\gamma)$ 
depend only on the endpoints of the oriented path $\gamma$ in
$\Gamma(X)$, if $G$ is simply-laced (that is, if all roots have the
same length for an appropriate choice of the $W$-invariant scalar
product on $\cX$; equivalently, $\Phi$ is a product of simple root
systems of type $A$, $D$ or $E$.) This assumption cannot be omitted,
as shown by
\medskip

\noindent
{\sl Example 3.} Let $G={\rm SP}(4)$ be the subgroup of 
${\rm GL}(4)$ preserving a non-degenerate symplectic form,
and let $H={\rm GL}(2)$ be the subgroup of $G$ preserving two
complementary lagrangian planes. The normalizer $N_G(H)$ contains $H$ 
as a subgroup of index 2. The graph $\Gamma(G/H)$ is as follows:

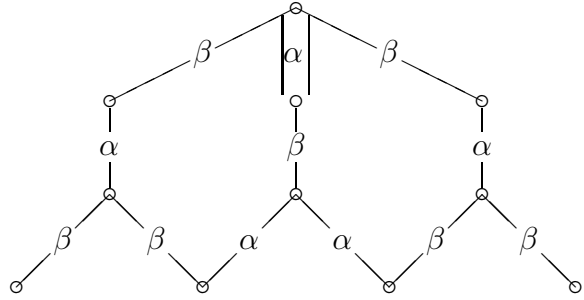
\begin{diagram}[abut,height=1.5em,width=1.5em]
&&&&&&\circ&&&&&\\
&&&&&\ldLine(4,2)~\beta&\dLine\alpha\dLine&\rdLine(4,2)~\beta&&&&\\
&&\circ&&&&\circ&&&&\circ&&\\
&&\dLine~\alpha&&&&\dLine~\beta&&&&\dLine~\alpha&\\
&&\circ&&&&\circ&&&&\circ&&\\
&\ldLine~\beta&&\rdLine~\beta&&\ldLine~\alpha&&\rdLine~\alpha&&
\ldLine~\beta&&\rdLine~\beta\\
\circ&&&&\circ&&&&\circ&&&&\circ\\
\end{diagram}

And here is $\Gamma(G/N_G(H))$:

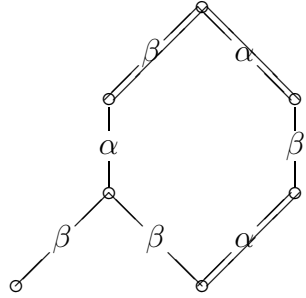
\begin{diagram}[abut,height=1.5em,width=1.5em]
&&&&\circ&&&&\\
&&&\ldLine~\beta\;\ldLine&&\rdLine~\alpha\;\rdLine\\
&&\circ&&&&\circ&\\
&&\dLine~\alpha&&&&\dLine~\beta\\
&&\circ&&&&\circ\\
&\ldLine~\beta&&\rdLine~\beta&&\ldLine~\alpha\;\ldLine\\
\circ&&&&\circ&&\\
\end{diagram}

\noindent
Using parabolic induction, one constructs similar examples for $\Phi$
of type $B$, $C$ or $F$.

\medskip

To proceed, we need the following definition taken from \cite{B1}:

\medskip

\noindent
{\sl Definition.} For $Y\in\cB(X)$, let $W(Y)$ be the set of all 
$w\in W$ such that the morphism 
$\pi_{Y,w}:\overline{BwB}\times^B Y\to GY$ is surjective
and generically finite. For $w\in W(Y)$, let $d(Y,w)$ be the degree of
$\pi_{Y,w}$.
\medskip

In other words, $W(Y)$ consists of all $w(\gamma)$ where $\gamma$ is
an oriented path from $Y$ to $GY$; moreover,
$d(Y,w(\gamma))=2^{\ell_N(\gamma)}$. By Lemma \ref{length},
$w^{-1}\in W^{\Delta(X)}$ for all $w\in W(Y)$.
\medskip

We now introduce a notion of neighbors in $W(Y)$, and we show that any
two elements of that set are connected by a chain of neighbors.
Let $\alpha$, $\beta$ be distinct simple roots and let $m$ be a
positive integer. Let
$$
(s_{\alpha}s_{\beta})^{(m)}=\cdots
s_{\beta}s_{\alpha}s_{\beta}s_{\alpha}
\eqno(m~{\rm terms}.)
$$
Then we have the braid relation 
$(s_{\alpha}s_{\beta})^{(m(\alpha,\beta))}=
(s_{\beta}s_{\alpha})^{(m(\alpha,\beta))}$,
where $m(\alpha,\beta)$ denotes the order of $s_{\alpha}s_{\beta}$ in
$W$.
\medskip

\noindent
{\sl Definition.} Two elements $u$ and $v$ of $W$ are
{\sl neighbors} if there exist $x$, $y$ in $W$ together with distinct
$\alpha$, $\beta$ in $\Delta$ and a positive integer
$m<m(\alpha,\beta)$ such that 
$$
u=x(s_{\alpha}s_{\beta})^{(m)}y,~
v=x(s_{\beta}s_{\alpha})^{(m)}y,~
{\rm and}~\ell(u)=\ell(x)+m+\ell(y)=\ell(v).
$$

For example, any two simple reflections are neighbors.

\begin{proposition}\label{neigh} 
Let $Y\in\cB(X)$ and let $u$, $v$ be
distinct elements of $W(Y)$. Then there exists a sequence
$(u=u_0,u_1,\ldots,u_n=v)$ in $W(Y)$ such that each $u_{i+1}$ is a
neighbor of $u_i$. 
\end{proposition}

\begin{proof} By induction on $\ell(u)=\ell(v)=\ell$, the case where
$\ell=1$ being evident. If there exists $\alpha\in\Delta$ such that
$\ell(us_{\alpha})=\ell(vs_{\alpha})=\ell-1$, then $P_{\alpha}$ raises
$Y$, and $us_{\alpha}$, $vs_{\alpha}$ are in $W(P_{\alpha}Y)$. Now the
induction assumption for $P_{\alpha}Y$ concludes the proof in this case.
Otherwise, we can find distinct $\alpha,\beta\in\Delta$ such that
$\ell(us_{\alpha})=\ell(vs_{\beta})=\ell-1$. Then $P_{\alpha}$ and
$P_{\beta}$ raise $Y$ to subvarieties of $P_{\alpha,\beta}Y$. Let $m$
be the common codimension of $P_{\alpha}Y$ and $P_{\beta}Y$ in
$P_{\alpha,\beta}Y$, then we have
$$
P_{\alpha,\beta}Y=\cdots P_{\alpha}P_{\beta}P_{\alpha}Y=
\overline{B\cdots s_{\alpha}s_{\beta}s_{\alpha}Y}\eqno (m~{\rm terms})
$$
Choose $x\in W(P_{\alpha,\beta}Y)$, then $W(Y)$ contains
$x(s_{\alpha}s_{\beta})^{(m)}$ and, similarly,
$x(s_{\beta}s_{\alpha})^{(m)}$, as neighbors. Moreover,
$W(P_{\alpha}Y)$ contains $us_{\alpha}$ and 
$x(s_{\beta}s_{\alpha})^{(m-1)}$, whereas $W(P_{\beta}Y)$ contains
$x(s_{\beta}s_{\alpha})^{(m-1)}$ and $vs_{\beta}$. Now we conclude by
the induction assumption for $P_{\alpha}Y$ and $P_{\beta}Y$.
\end{proof}

Neighbors in $W(Y)$ are also close to each other for the
Bruhat-Chevalley order $\leq$ on $W$:

\begin{proposition}\label{sup}
Let $Y\in\cB(X)$. For any neighbors $u,v\in W(Y)$, there exists 
$w\in W$ such that $u\leq w$, $v\leq w$, $w^{-1}\in W^{\Delta(X)}$ and
$\ell(w)=\ell(u)+1=\ell(v)+1$.
\end{proposition}

\begin{proof}
Write $u=x(s_{\alpha}s_{\beta})^{(m)}y$ and
$v=x(s_{\beta}s_{\alpha})^{(m)}y$. Let
$$
w=x(s_{\alpha}s_{\beta})^{(m)}s_{\beta}y.
$$
We claim that $\ell(w)$ equals
$\ell(x)+m+1+\ell(y)=\ell(u)+1=\ell(v)+1$. Otherwise,
$\ell(w)\leq \ell(x)+\ell(y)+m-1<l(u)$ and
$w=uy^{-1}s_{\beta}y=us_{y^{-1}(\beta)}$. By the strong exchange
condition (\cite{H} Theorem 5.8 applied to $u$), one of the following
cases occurs: 

\noindent
(i) $w=x'(s_{\alpha}s_{\beta})^{(m)}y$ where
$\ell(x')=\ell(x)-1$. Comparing both expressions for $w$, we obtain 
$x'(s_{\alpha}s_{\beta})^{(m)}=x(s_{\alpha}s_{\beta})^{(m)}s_{\beta}$.
Thus, there exists $\gamma\in\Phi^+_{\alpha,\beta}$ such that
$x'=xs_{\gamma}$. But $\ell(xs_{\alpha})=\ell(xs_{\beta})=\ell(x)+1$,
for $\ell(x(s_{\alpha}s_{\beta})^{(m)}y)=
\ell(x(s_{\beta}s_{\alpha})^{(m)}y)=\ell(x)+m+\ell(y)$. It follows
that $x(\alpha)$ and $x(\beta)$ are in $\Phi^+$. Thus, 
$x\in W^{\alpha,\beta}$. Since $s_{\gamma}\in W_{\alpha,\beta}$, we
have $\ell(x')=\ell(x)+\ell(s_{\gamma})\geq \ell(x)$, a
contradiction.

\noindent
(ii) $w=xzy$ where $z$ is obtained from
$(s_{\alpha}s_{\beta})^{(m)}$ by deleting a simple reflection. Then
the equality
$z=(s_{\alpha}s_{\beta})^{(m)}s_{\beta}$ leads to a braid
relation of length at most $m<m(\alpha,\beta)$, a contradiction.

\noindent
(iii) $w=x(s_{\alpha}s_{\beta})^{(m)}y'$ where
$\ell(y')=\ell(y)-1$. Then $y'=s_{\beta}y$. But
$\ell(s_{\beta}y)=\ell(y)+1$, for $\ell(v)=\ell(x)+m+\ell(y)$; a
contradiction.

By the claim and \cite{H} Theorem 5.10, we have $u\leq w$ and 
$v\leq w$. Write $w=w''w'$ where
$w''\in W_{\Delta(X)}$ and $(w')^{-1}\in W^{\Delta(X)}$; then
$\ell(w)=\ell(w')+\ell(w'')$. Since $u^{-1}\leq w^{-1}$ and 
$u^{-1}\in W^{\Delta(X)}$, it follows that $u^{-1}\leq(w')^{-1}$ by
\cite{Deo} Lemma 3.5. Thus, $u\leq w'$ and $v\leq w'$. Since $u\neq v$
and $\ell(u)=\ell(v)=\ell(w)-1\geq \ell(w')-1$, we must have $w=w'$,
so that $w^{-1}\in W^{\Delta(X)}$.
\end{proof}

Recall that $r(Y)\leq r(X)$ for any $Y\in\cB(X)$, see \cite{K2} 
Corollary 2.4. If equality holds, then neighbors in $W(Y)$ have 
a very simple form: 

\begin{proposition}\label{equal}
Let $Y\in\cB(X)$ such that $r(Y)=r(X)$; let $u,v\in W(Y)$ be
neighbors. Then $u=xs_{\alpha}y$ and $v=xs_{\beta}y$ where $x,y\in W$
and $\alpha,\beta$ are orthogonal simple roots such that
$\ell(u)=\ell(v)=\ell(x)+\ell(y)+1$. Moreover, $\cX(X)$ contains
$x(\alpha+\beta)$.
\end{proposition}

\begin{proof}
Write $u=x(s_{\alpha}s_{\beta})^{(m)}y$ and
$v=x(s_{\beta}s_{\alpha})^{(m)}y$ as in the definition of neighbors. 
Then $x(s_{\alpha}s_{\beta})^{(m)}$ and $x(s_{\beta}s_{\alpha})^{(m)}$ 
are neighbors in $W(\overline{ByY})$. Moreover, 
$r(\overline{ByY})\geq r(Y)$, whence $r(\overline{ByY})=r(X)$. Thus, 
we may assume that $y=1$.

Let $Y'=\overline{B(s_{\alpha}s_{\beta})^{(m)}Y}$
and $Y''=\overline{B(s_{\beta}s_{\alpha})^{(m)}Y}$, then we obtain
similarly: $r(Y')=r(Y'')=r(X)$ and $x\in W(Y')\cap W(Y'')$. If 
$x\neq 1$, write $x=s_{\gamma}x'$ where $\gamma\in\Delta$ and
$\ell(x)=\ell(x')+1$. Then $\overline{Bx'Y'}$ and $\overline{Bx'Y''}$
have rank $r(X)$ and are raised to $X$ by $\gamma$. Thus,
$\overline{Bx'Y'}=\overline{Bx'Y''}$ and, by induction on $\ell(x)$,
we obtain $Y'=Y''$. This subvariety is stable under
$P_{\alpha,\beta}$. Applying Lemmas \ref{reduc} and 
\ref{nor}, we may assume that $Y'=X$ (i.e., $x=1$),
$\Delta=\{\alpha,\beta\}$ and $X=G/H$ where the center of $G$ is
trivial and $H$ has finite index in its normalizer. Moreover, we have
$P(X)=B$, for $P_{\alpha}$ and $P_{\beta}$ do not stabilize $X^0$.

We claim that any $Z\in\cB(X)$ can be written as
$$
\overline{B(s_{\alpha}s_{\beta})^{(n)}Y}=\cdots P_{\beta}P_{\alpha}Y
~{\rm or}~
\overline{B(s_{\beta}s_{\alpha})^{(n)}Y}=\cdots P_{\alpha}P_{\beta}Y
\eqno(n~{\rm terms}),
$$
where $n=\dim(Z)-\dim(Y)$ satisfies $0\leq n\leq m$. For this, we argue 
by induction on the codimension of $Z$ in $X$. We may assume that 
$\alpha$ raises $Z$. By the induction assumption, we have
$$
P_{\alpha}Z=P_{\beta}P_{\alpha}\cdots Y~{\rm or}~
P_{\alpha}Z=P_{\alpha}P_{\beta}\cdots Y
\eqno(n+1~{\rm terms}).
$$ 
In the latter case, let $Z'=P_{\beta}\cdots Y$ ($n$ terms). Since 
$P_{\alpha}Z=P_{\alpha}Z'$ and $r(Z)=r(Z')=r(P_{\alpha}Z)=r(Y)$, 
it follows that $Z=Z'$. In the former case, $P_{\alpha}Z$ is stable 
under $G$ and hence equal to $X$;  in particular, $Z$ has codimension 
$1$ in $X$. Now $X=P_{\alpha}P_{\beta}\cdots Y$ ($m$ terms), so that 
we are in the previous case.

By the claim, all $B$-orbit closures in $X$ have the same rank, and
$Y^0$ is the unique closed $B$-orbit. Let $y\in Y^0$; we may assume
that $H=G_y$. Since the $H$-orbit in $G/B$ corresponding to the
$B$-orbit $Y^0$ in $G/H$ is closed, the connected isotropy group
$B_y^0$ is a Borel subgroup of $H^0$. It follows that
$r(Y)=r(B)-r(B_y)=2-r(H)$. On the other hand, $r(Y)=r(G/H)$ by
assumption. Thus, $r(G/H)=2-r(H)$.

If $r(G/H)=0$ then $H$ is a parabolic subgroup of $G$ (in fact, a
Borel subgroup as $P(G/H)=B$.) Moreover, $Y$ is the $B$-fixed point in
$G/H$. But then $W(Y)$ consists of a unique element (of maximal length
in $W$), a contradiction.

If $r(G/H)=1$ then $r(H)=1$ as well. Using the classification of
homogeneous spaces of rank $1$ under semi-simple groups of rank $2$
(see e.g. Table 1 of \cite{W}), this forces 
$G={\rm PGL}(2)\times{\rm PGL}(2)$ and $H={\rm PGL}(2)$ embedded
diagonally in $G$. As a consequence, the simple roots $\alpha$ and
$\beta$ are orthogonal, and $\cX(G/H)$ is generated by
$\alpha+\beta$.

If $r(G/H)=2$ then $r(H)=0$, that is, $H^0$ is unipotent. Since $G/H$
is spherical, $H^0$ is a maximal unipotent subgroup of $G$. This
contradicts the assumption that $H$ has finite index in its
normalizer.
\end{proof}

\begin{proposition}\label{constant}
If $G$ is simply-laced, then

\noindent
(i) for any oriented path $\gamma$ in $\Gamma(X)$, both
$\ell_T(\gamma)$ and $\ell_N(\gamma)$ depend only on the endpoints of
$\gamma$.

\noindent
(ii) for any $Y\in\cB(X)$, there exists an oriented path $\gamma$
joining $Y$ to $X$ through a sequence of simple edges followed by a
sequence of double edges.
\end{proposition}

\begin{proof}
(i) Let $Y$ (resp. $Y'$) be the source (resp. target) of $\gamma$, and
let $\delta$ be another oriented path from $Y$ to $Y'$. By Lemma
\ref{length}, it suffices to show that
$\ell_N(\gamma)=\ell_N(\delta)$. Joining $Y'$ to $X$ by an oriented
path, we reduce to the case where $Y'=X$; then 
$w(\gamma)$ and $w(\delta)$ are in $W(Y)$. By Proposition \ref{neigh},
we may assume moreover that $w(\gamma)$ and $w(\delta)$ are
neighbors. Using Lemmas \ref{reduc} and \ref{nor}, we 
reduce to the case where the center of $G$ is trivial, 
$\Delta=\{\alpha,\beta\}$, $X=G/H$ where $H$ has finite index in its
normalizer, $w(\gamma)=(s_{\alpha}s_{\beta})^{(m)}$ and
$w(\delta)=(s_{\beta}s_{\alpha})^{(m)}$ for some $m<m(\alpha,\beta)$. 

Since $G$ is simply-laced, we have either 
$G={\rm PGL}(2)\times{\rm PGL}(2)$ and $m(\alpha,\beta)=2$, or 
$G={\rm PGL}(3)$ and $m(\alpha,\beta)=3$. In particular, $m\leq 2$. 
If $m=1$ then $\ell_N(\gamma)=\ell_N(\delta)=0$ by Proposition
\ref{type}. If $m=2$ then $G={\rm PGL}(3)$. Using Lemma \ref{length}
(iv), we may assume moreover that $H$ is not contained in any Borel
subgroup. Then we see by inspection that $H$ is conjugate to 
${\rm PO}(3)$ or to ${\rm GL}(2)$.

In the latter case, here is $\Gamma(G/H)$:

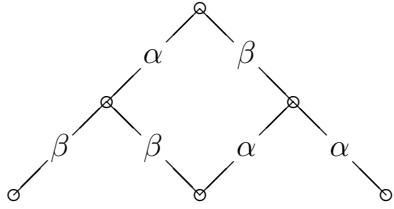
\begin{diagram}[abut,height=1.5em,width=1.5em]
&&&&\circ&&&&\\
&&&\ldLine~\alpha&&\rdLine~\beta&&\\
&&\circ&&&&\circ&\\
&\ldLine~\beta&&\rdLine~\beta&&\ldLine~\alpha&&\rdLine~\alpha\\
\circ&&&&\circ&&&&\circ\\
\end{diagram}

Thus, $\ell_N(\gamma)=\ell_N(\delta)=0$. 

In the former case, we have $\ell_N(\gamma)=\ell_N(\delta)=1$, 
since $\Gamma(G/H)$ is as follows:

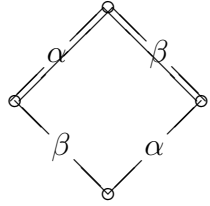
\begin{diagram}[abut,height=1.5em,width=1.5em]
  &&\circ &&\\
&\ldLine~\alpha\;\ldLine&&\rdLine\;\rdLine~\beta&\\
\circ &&&& \circ \\
& \rdLine~\beta && \ldLine~\alpha &\\
&& \circ &&
\end{diagram}

(ii) Let $\gamma$ be an oriented path joining $Y$ to $X$. We may
assume that $\gamma$ contains double edges. Consider the lowest
maximal subpath $\delta$ of $\gamma$ that consists of double edges
only; we may assume that the endpoint of $\delta$ is not $X$. Let $Y'$
be the source of the top edge of $\delta$, and let $\alpha$
(resp. $\beta$) be the label of that edge (resp. of the next edge of
$\gamma$, a simple edge by assumption.) We claim that there exists an
oriented path $\gamma'$ joining $Y'$ to $X$ and beginning with a
simple edge; then assertion (ii) will follow by induction on
$\ell(\delta)+{\rm codim}_X(Y')$.

To check the claim, it suffices to join $Y'$ to $P_{\alpha\beta}Y'$ by
an oriented path $\gamma$' beginning with a simple edge. As above, we 
reduce to the case where $G$ equals ${\rm PGL}(2)\times{\rm PGL}(2)$
or ${\rm PGL}(3)$, and $H$ is not contained in a Borel subgroup of
$G$; Moreover, $H$ has finite index in its normalizer. Using the fact
that $\Gamma(G/H)$ contains a double edge followed by a simple edge,
one checks that $H$ is a product of subgroups of ${\rm PGL}(2)$ if
$G={\rm PGL}(2)\times{\rm PGL}(2)$; and if $G={\rm PGL}(3)$, then $H$
is conjugate to the subgroup of Example 1, or to its transpose. The
path $\gamma'$ exists in all these cases.
\end{proof}

From Proposition \ref{constant} we will deduce a criterion for the
graph of a spherical variety to contain simple edges only. To
formulate it, we need more notation, and a preliminary result.

Let $D\in\cD(X)$ be a color; then $D$ is the closure of its
intersection with the open $G$-orbit $G/H$. Let $\tilde D$ be the
preimage in $G$ of $D\cap G/H$. Replacing $G$ by a finite cover, we
may assume that $\tilde D$ is the divisor of a regular function $f_D$
on $G$. Then $f_D$ is an eigenvector of $B$ acting by left
multiplication; let $\omega_D$ be its weight. Since $f_D$ is uniquely
defined up to multiplication by a regular invertible function on $G$,
then $\omega_D$ is unique up to addition of a character of $G$. 
In particular, for any $\alpha\in\Delta$, the number 
$\langle\omega_D,\check\alpha\rangle$
is a non-negative integer depending only on $D$ and $\alpha$.

\begin{lemma}\label{degree}
(i) The degree $d(D,\alpha)$ of the morphism
$\pi_{D,\alpha}:P_{\alpha}\times^B D\to X$ equals
$\langle\omega_D,\check\alpha\rangle$ if $\pi_{D,\alpha}$ is
generically finite; otherwise, $\langle\omega_D,\check\alpha\rangle=0$.
\noindent

(ii) For any $G$-orbit closure $X'$ in $X$ and for any 
$D'\in\cD(X')$, there exists $D\in\cD(X)$ such that $D'$ is an 
irreducible component of $D\cap X'$. Then
$\langle\omega_{D'},\check\alpha\rangle\leq
\langle\omega_D,\check\alpha\rangle$
for all $\alpha\in\Delta$.
\end{lemma}

\begin{proof}
(i) Note that $D$ is $P_{\alpha}$-stable if and only if $f_D$ is an 
eigenvector of $P_{\alpha}$, that is, $\omega_D$ extends to a 
character of that group. This amounts to:
$\langle\omega_D,\check\alpha\rangle=0$.

Let $V$ be the $H$-stable divisor in $G/B$ corresponding to the
$B$-stable divisor $D\cap G/H$. Then
$V$ is the zero scheme of a section of the homogeneous line bundle on
$G/B$ associated with the character $\omega_D$ of $B$. Let 
$p:G/B\to G/P_{\alpha}$ be the natural map, then $d(D,\alpha)$ equals 
the degree of the restriction $p_V:V\to G/P_{\alpha}$. The latter 
degree is the intersection number of $V$ with a fiber of $p$, that is,
$\langle\omega_D,\check\alpha\rangle$.

(ii) For the first assertion, it suffices to show existence of 
$D\in\cD(X)$ containing $D'$ and not containing $X'$; 
but this follows from \cite{K0} Theorem 3.1. For the second assertion,
note that $P_{\alpha}$ stabilizes $D'$ if it stabilizes $D$. Thus, 
$\langle\omega_{D'},\check\alpha\rangle=0$ if 
$\langle\omega_D,\check\alpha\rangle=0$. On the other hand, if 
$\langle\omega_D,\check\alpha\rangle=1$ then $\pi_{D,\alpha}$ is
birational. Restricting to $P_{\alpha}\times^B D'$, it follows that
$\pi_{D',\alpha}$ is birational if generically finite.
\end{proof}

A direct consequence of Lemma \ref{degree} and Proposition
\ref{constant} is 

\begin{corollary}\label{simple}
If $G$ is simply-laced, then the following conditions are equivalent:

\noindent
(i) Each edge of $\Gamma(X)$ is simple.

\noindent
(ii) For any $D\in\cD(X)$ and $\alpha\in\Delta$, we have
$\langle\omega_D,\check\alpha\rangle \leq 1$.
\end{corollary}

This criterion applies, e.g., to all embeddings of the following
symmetric spaces:  
${\rm GL}(p+q)/{\rm GL}(p)\times{\rm GL}(q)$, 
${\rm SL}(2n)/{\rm SP}(2n)$, ${\rm SO}(2n)/{\rm GL}(n)$ and $E_6/F_4$.
For this, one uses the explicit description of colors of symmetric
spaces given in \cite{V}. Further applications will be given after
Theorem \ref{normal} below. 

Note that Corollary \ref{simple} does not extend to multiply-laced 
groups $G$. Consider, for example, $G={\rm SO}(2n+1)$ and its subgroup 
$H={\rm O}(2n)$, the stabilizer of a non-degenerate line in $\mC^{2n+1}$. 
Then the homogeneous space $G/H$ is spherical of rank $1$ and its graph
consists of a unique oriented path: a double edge followed by $n-1$
simple edges.

\section{Orbit closures in regular varieties}

Recall from \cite{BDP} that a variety $X$ with an action of $G$ is
called {\sl regular} if it satisfies the following three conditions:  

\noindent
(i) $X$ is smooth and contains a dense $G$-orbit whose complement is a
union of irreducible smooth divisors (the {\sl boundary divisors})
with normal crossings.

\noindent
(ii) Any $G$-orbit closure in $X$ is the transversal intersection of
those boundary divisors that contain it.

\noindent
(iii) For any $x\in X$, the normal space to the orbit $Gx$ contains a
dense orbit of the isotropy group of $x$.
\medskip

Any regular $G$-variety $X$ contains only finitely many $G$-orbits.
Their closures are the $G$-stable subvarieties of $X$; they are
regular $G$-varieties as well.

Regular varieties are closely related with spherical varieties: any
complete regular $G$-variety is spherical, and any spherical
$G$-homogeneous space $G/H$ admits an open equivariant embedding into
a complete regular $G$-variety $X$, see \cite{BB} 2.2.

Let $Z$ be a closed $G$-orbit in complete regular $X$, then the
isotropy group of each point of $Z$ is a parabolic subgroup of
$G$. Thus, $Z$ contains a unique $T$-fixed point $z$ such that $Bz$ is
open in $Z$; we shall call $z$ the base point of $Z$. In fact, the
isotropy group $Q=G_z$ is opposed to $P(X)$, see e.g. \cite{BB} 2.2.
 
We next recall the local structure of complete regular varieties, see 
e.g. \cite{BB} 2.3. For such a variety $X$, set $P=P(X)$ and
$L=L(X)$. Let $X_0$ be the set of all $x\in X$ such that $Bx$ is open
in $Gx$. Then $X_0$ is an open $P$-stable subset of $X$: the
complement of the union of all colors. Moreover, there exists an
$L$-stable subvariety $S$ of $X_0$, fixed pointwise by $[L,L]$, such
that the map
$$
\begin{matrix}
R_u(P)\times S&\to& X_0\cr
(g,x)&\mapsto&gx\cr
\end{matrix}
$$
is an isomorphism. As a consequence, $S$ is a smooth toric variety 
(for a quotient of $T$) of dimension $r(X)$, the rank of $X$;
moreover, $S$ meets each $G$-orbit along a unique $T$-orbit.
Let $\varphi:X_0\cong R_u(P)\times S\to S$ be the second projection,
then $\varphi$ is $L$-equivariant; it can be seen as the quotient map
by the action of $R_u(P)$.

We now turn to $B$-orbit closures. Let $Y\in\cB(X)$; since $GY$ is
regular, we may assume that $GY=X$. Then, by \cite{B1} 1.4, $Y$ meets
all $G$-orbit closures properly; moreover, for any closed $G$-orbit
$Z$, the irreducible components of $Y\cap Z$ are the Schubert
varieties $\overline{Bw^{-1}z}$ where $w\in W(Y)$, and the
intersection multiplicity of $Y$ and $Z$ along $\overline{Bw^{-1}z}$
equals $d(Y,w)$. To describe the intersection of $Y$ with
arbitrary $G$-orbit closures, we shall study the local structure of
$Y$ along $\overline{Bw^{-1}z}$ for a fixed $w\in W(Y)$. It will be
more convenient to consider the translate $wY$ along
$\overline{wBw^{-1}z}$. 

Note that $wY$ meets $X_0$ (because $\overline{BwY}=X$), and that the
intersection $wY\cap X_0$ is stable by the group $wBw^{-1}\cap P$. The
latter contains $R_u(P)\cap wUw^{-1}$ as a normal subgroup. We shall
see that $R_u(P)\cap wUw^{-1}$ acts freely on $wY\cap X_0$, with
section
$$
S_{Y,w}=wY\cap (U\cap wU^-w^{-1})S.
$$
Note that $U\cap wU^-w^{-1}$ is contained in $R_u(P)$, because
$w^{-1}\in W^P$. Thus, $S_{Y,w}$ is a closed $T$-stable subvariety of
$wY\cap X_0$. Let 
$$
\varphi_{Y,w}:S_{Y,w}\to S
$$
be the restriction of $\varphi:X_0\to S$, then $\varphi_{Y,w}$ is
$T$-equivariant.

\begin{proposition}\label{slice} 
Keep notation as above.

\noindent
(i) The map
$$\begin{matrix}
(R_u(P)\cap wUw^{-1})\times S_{Y,w}&\to&wY\cap X_0\cr
(g,x)&\mapsto& gx\cr
\end{matrix}$$
is an isomorphism.

\noindent
(ii) The variety $S_{Y,w}$ is irreducible and meets each $G$-orbit
along a unique $T$-orbit. In particular, $S_{Y,w}\cap GY^0$ is a
unique $T$-orbit, dense in $S_{Y,w}$ and contained in $wY^0$; and
$S_{Y,w}\cap Z=\{z\}$ for any closed $G$-orbit $Z$ with base
point $z$.

\noindent
(iii) The morphism $\varphi_{Y,w}$ is finite surjective of degree
$d(Y,w)$.
\end{proposition}

\begin{proof}
(i) The product map 
$(R_u(P)\cap wUw^{-1})\times (R_u(P)\cap wU^-w^{-1})\to R_u(P)$
is an isomorphism; moreover, 
$R_u(P)\cap wU^-w^{-1}=U\cap wU^-w^{-1}$. Therefore, the product map 
$$
(R_u(P)\cap wUw^{-1})\times (U\cap wU^-w^{-1})S\to X_0
$$
is an isomorphism. The assertion follows by intersecting with $wY$.

(ii) and (iii) The union of all $G$-orbits in $X$ that contain $Z$ in
their closure is a $G$-stable open subset of $X$. Thus, we may assume
that $Z$ is the unique closed $G$-orbit in $X$. Let 
$D_1,\ldots,D_r$ be the boundary divisors, then $r=r(X)$. Moreover,
$S$ is isomorphic to affine space $\mA^r$ with coordinate functions
$x_1,\ldots,x_r$, equations of $D_1\cap S,\ldots,D_r\cap S$. The
compositions $f_1=x_1\circ\varphi,\ldots,f_r=x_r\circ\varphi$ are
equations of $D_1\cap X_0,\ldots,D_r\cap X_0$; they generate the
ideal of $Z\cap X_0=Bz$ in $X_0$. The map $\varphi:X_0\to S$
identifies to $(f_1,\ldots,f_r):X_0\to\mA^r$. The intersections of
$G$-orbit closures with $X_0$ are the pull-backs of coordinate
subspaces of $\mA^r$.

By (i), $S_{Y,w}$ is irreducible. We check that 
$S_{Y,w}\cap Z=\{z\}$. For this, note that the product map
$$
(R_u(P)\times wUw^{-1})\times (S_{Y,w}\cap Z)\to 
wY\cap X_0\cap Z = wY\cap Bz
$$
is an isomorphism. Moreover, since $Y$ meets $Z$ properly, with
$\overline{Bw^{-1}z}$ as an irreducible component, it follows that 
$wY\cap Bz$ is equidimensional, with 
$\overline{wBw^{-1}z}\cap Bz=(B\cap wBw^{-1})z$
as an irreducible component. The latter is isomorphic to $R_u(P)\cap
wUw^{-1}$. Thus, the $T$-stable set $S_{Y,w}\cap Z$ is finite, so that
it consists of $T$-fixed points. Since $z$ is the unique $T$-fixed
point in $Bz$, our assertion follows. 

The map $\varphi_{Y,w}:S_{Y,w}\to S$ identifies with
$(f_1,\ldots,f_r):S_{Y,w}\to\mA^r$. We just saw that the
set-theoretical fiber of $0$ is $\{z\}$. Since $0$ is the unique
closed $T$-orbit in $\mA^r$, all fibers of $\varphi_{Y,w}$ are
finite. Thus, $S_{Y,w}$ contains a dense $T$-orbit. Since $S_{Y,w}$ is
affine and contains a $T$-fixed point $z$, it follows that
$\varphi_{Y,w}$ is finite and that the pull-back of any $T$-orbit in
$S$ is a unique $T$-orbit. This implies (ii).

Finally, we check that the degree of $\varphi_{Y,w}$ equals
$d(Y,w)$, that is, the degree of the natural map 
$\overline{BwB}\times ^B Y\to X$. For this, note that the map 
$$
U\cap wU^-w^{-1}\to\overline{BwB}/B,~g\mapsto gwB/B
$$ 
is an open immersion. Thus, $d(Y,w)$ is the degree of the product map
$(U\cap wU^-w^{-1})\times wY\to X$, or, equivalently, of its restriction
$$
p:(U\cap wU^-w^{-1})\times (wY\cap X_0)\to X_0.
$$
The latter map fits into a commutative diagram
$$\begin{matrix}
(U\cap wU^-w^{-1})\times(wY\cap X_0)&\to&X_0\cr
\downarrow&&\downarrow\cr
S_{Y,w}&\to&S,\cr
\end{matrix}$$
where the bottom horizontal map is $\varphi_{Y,w}$ ; indeed, 
$$
(U\cap wU^-w^{-1})\times(wY\cap X_0)\cong
(R_u(P)\cap wU^-w^{-1})\times (R_u(P)\cap wUw^{-1})\times S_{Y,w}
$$
by (i). Moreover, the fibers of the right (resp. left) vertical map
are isomorphic to $R_u(P)$ (resp. to 
$(R_u(P)\cap wU^-w^{-1})\times (R_u(P)\cap wUw^{-1})\cong R_u(P)$.)
Thus, the diagram is cartesian, and the degree of $p$ equals the
degree of $\varphi_{Y,w}$.
\end{proof}

Thus, we can view $S_{Y,w}$ as a ``slice'' in $wY$ to 
$wBw^{-1}z=(R_u(P)\cap wUw^{-1})z$ at $z$. But $S_{Y,w}$ may be
non transversal to $wY$ at $z$: indeed, the intersection multiplicity
of $S_{Y,w}$ and $wY$ at $z$ equals the intersection multiplicity of
$Z$ and $Y$ along $\overline{Bw^{-1}z}$, and the latter equals
$d(Y,w)$ by \cite{B1} 1.4 (alternatively, this can be deduced from
Proposition \ref{slice} (iii).) On the other hand, it is not clear
whether $S_{Y,w}$ is smooth, that is, $Y\cap w^{-1}X_0$ consists of
smooth points of $Y$; see Corollary \ref{components} below for a
partial answer to this question.

We now relate the ``slices'' associated with both endpoints of an edge
in $\Gamma(X)$. Let $Y\in\cB(X)$ and let $\alpha\in\Delta$ raising
$Y$. Choose $v\in W(P_{\alpha}Y)$, then $w=vs_{\alpha}$ is in $W(Y)$,
and $\ell(w)=\ell(v)+1$. Thus, $v(\alpha)\in\Phi^+\cap w(\Phi^-)$. Let
$U_{v(\alpha)}$ be the corresponding unipotent subgroup of dimension
$1$, then $U_{v(\alpha)}$ is contained in $R_u(P)\cap vUv^{-1}$.

\begin{proposition}\label{step}
With notation as above, $S_{Y,w}$ is contained in
$U_{v(\alpha)}S_{P_{\alpha}Y,v}$, and the latter is isomorphic to
$U_{v(\alpha)}\times S_{P_{\alpha}Y,\tau}$. Denoting by 
$$
\varphi_{Y,\alpha}:S_{Y,w}\to S_{P_{\alpha}Y,v}
$$
the corresponding projection, then 
$\varphi_{Y,w}=\varphi_{P_{\alpha}Y,v}\circ \varphi_{Y,\alpha}$. 
Moreover, $\varphi_{Y,\alpha}$ is finite surjective of degree
$d(Y,\alpha)$.
\end{proposition}

\begin{proof}
We have
$$\displaylines{
S_{Y,w} = wY\cap (U\cap wU^-w^{-1})S 
= wY\cap U_{v(\alpha)}(U\cap vU^-v^{-1})S
\hfill\cr\hfill
\subseteq vP_{\alpha}Y\cap U_{v(\alpha)}(U\cap vU^-v^{-1})S
=U_{v(\alpha)}(vP_{\alpha}Y\cap (U\cap vU^-v^{-1})S)
=U_{v(\alpha)}S_{P_{\alpha}Y,v}.
\cr}$$
Moreover, since $U_{v(\alpha)}\subseteq R_u(P)\cap vUv^{-1}$, the
product map 
$U_{v(\alpha)}\times S_{P_{\alpha}Y,v}\to
U_{v(\alpha)}S_{P_{\alpha}Y,v}$
is an isomorphism. Now the equality 
$\varphi_{Y,w}=\varphi_{P_{\alpha}Y,v}\circ \varphi_{Y,\alpha}$
follows from the definitions. Together with Proposition \ref{slice}
(iii), it implies that $\varphi_{Y,\alpha}$ is
finite surjective of degree $d(Y,w)d(P_{\alpha}Y,v)^{-1}=d(Y,\alpha)$.
\end{proof}

Using Proposition \ref{slice}, we analyze the intersection
of $Y$ with an arbitrary $G$-orbit closure, generalizing \cite{B1}
Theorem 1.4.

\begin{theorem}\label{indices}
Let $X$ be a complete regular $G$-variety, let $Y\in\cB(X)$ be such
that $GY=X$ and let $X'$ be a $G$-orbit closure. Then $W(Y)$ is the
disjoint union of the $W(C)$ where $C$ runs over all irreducible
components of $Y\cap X'$. Moreover, for any such $C$ and $w\in W(C)$, 
we have 
$$
d(Y,w)=d(C,w)\;i(C,Y\cdot X';X)
$$
where $i(C,Y\cdot X';X)$ denotes the intersection multiplicity of $Y$
and $X'$ along $C$ in $X$. As a consequence, this multiplicity is a
power of $2$.
\end{theorem}

\begin{proof}
By \cite{B1} Lemma 1.3, $W(Y)$ is the union of the $W(C)$. Choose
$C$ and $w\in W(C)$, then $C\cap w^{-1}X_0$ is an irreducible
component of $Y\cap w^{-1}X_0\cap X'$. The latter is isomorphic to
$(U\cap w^{-1}R_u(P))\times w^{-1}(S_{Y,w}\cap X')$, and 
$S_{Y,w}\cap X'$ is a unique $T$-orbit, by Proposition
\ref{slice}. It follows that $Y\cap w^{-1}X_0\cap X'=C\cap w^{-1}X_0$ 
is irreducible, so that $C$ is uniquely determined by $w$. Equivalently,
the $W(C)$ are pairwise disjoint.

Let $Z$ be a closed $G$-orbit in $X'$, then 
$$
d(Y,w)=i(\overline{Bw^{-1}z},Y\cdot Z;X)
=i(\overline{Bw^{-1}z}\cap w^{-1}X_0,
(Y\cap w^{-1}X_0)\cdot(Z\cap  w^{-1}X_0); w^{-1}X_0),
$$
where the former equality follows from \cite{B1} 1.4, and the latter
from \cite{F} 8.2. Moreover, we have by Proposition \ref{slice}: 
$\overline{Bw^{-1}z}\cap w^{-1}X_0=Bw^{-1}z$ and 
$Z\cap w^{-1}X_0=w^{-1}Bz$. Thus,
$$
d(Y,w)=i(Bw^{-1}z,(Y\cap w^{-1}X_0)\cdot w^{-1}Bz,w^{-1}X_0).
$$
Using the fact that
$Y\cap w^{-1}X_0\cap X'=C\cap w^{-1}X_0$ is irreducible, together with
associativity of intersection multiplicities (see \cite{F} 7.1.8),
we obtain
$$\displaylines{
d(Y,w)=i(Bw^{-1}z,(C\cap w^{-1}X_0)\cdot w^{-1}Bz;w^{-1}X_0\cap X')
\;i(C,Y\cdot X';X)
\hfill\cr\hfill
=i(\overline{Bw^{-1}z},C\cdot Z;X')\;i(C,Y\cdot X';X)
=d(C,w)\;i(C,Y\cdot X';X).
\cr}$$
\end{proof}

These results motivate the following
\medskip

\noindent
{\sl Definition}. A $B$-orbit closure $Y$ in an arbitrary spherical
variety $X$ is {\sl multiplicity-free} if $d(Y,w)=1$ for all 
$w\in W(Y)$. Equivalently, the edges of all oriented paths in
$\Gamma(X)$ with source $Y$ are simple.
\medskip

For example, $Y$ is multiplicity-free if $r(Y)=r(GY)$, or if the
isotropy group in $G$ of a point of $Y^0$ is contained in a Borel
subgroup of $G$ (this follows from Lemma \ref{length}.) 

Other examples of multiplicity-free orbit closures arise from
parabolic induction: if $X=G\times^{P_I} X'$ is induced from $X'$ and
if $Y=\overline{BwY'}$ with $w\in W^I$ and $Y'\in\cB(X')$, then $Y$ is
multiplicity-free if and only if $Y'$ is (this follows from Lemma
\ref{induced} or, alternatively, from \cite{B1} 1.2).
\medskip

\begin{corollary}\label{components}
Let $X$ be a complete regular $G$-variety, $Y$ a multiplicity-free
$B$-stable subvariety such that $GY=X$, and $X'$ a $G$-orbit closure
in $X$. Then all irreducible components of $Y\cap X'$ are
multiplicity-free $B$-orbit closures of $X'$, and the corresponding
intersection multiplicities equal $1$. Moreover, for any $w\in W(Y)$,
the map $\varphi_{Y,w}:S_{Y,w}\to S$ is an isomorphism. As a
consequence, $Y\cap w^{-1}X_0$ consists of smooth points of $Y$.
\end{corollary}

\begin{proof}
The first assertion follows from Theorem \ref{indices}. By
Proposition \ref{slice}, $\varphi_{Y,w}$ is finite surjective of
degree $1$, hence an isomorphism because $S$ is smooth.
\end{proof}

Returning to arbitrary $B$-orbit closures in a complete regular
$G$-variety, we now show that their intersections with $G$-orbit
closures satisfy Hartshorne's connectedness theorem, see \cite{E}
18.2. That theorem is proved there for schemes of depth at least $2$;
but $B$-orbit closures may have depth $1$ at some points, see Example
4 in the next section.

\begin{theorem}\label{connected}
Let $X$ be a complete regular $G$-variety, $Y$ a $B$-orbit closure,
and $X'$ a $G$-orbit closure in $X$. Then $Y\cap X'$ is connected in
codimension $1$ (that is, the complement in $Y\cap X'$ of any closed
subset of codimension at least $2$ is connected.)
\end{theorem}

\begin{proof}
We may assume that $GY=X$. If $X'=Z$ is a closed $G$-orbit, then the
assertion follows from the description of $Y\cap Z$ in terms of
$W(Y)$, together with Propositions \ref{neigh} and \ref{sup}. Indeed,
for any $w\in W$ such that $w^{-1}\in W^{\Delta(X)}$, we have 
$\ell(w)=\ell(w^{-1})={\rm codim}_Z(\overline{Bw^{-1}z})$,
where $z$ is the base point of $Z$.

For arbitrary $X'$, let $Z$ be a closed $G$-orbit in $X'$. Let
$Y'_1$, $Y'_2$ be unions of irreducible components of $Y\cap X'$ such
that $Y\cap X'=Y'_1\cup Y'_2$. Then $Y'_1\cap Z$ and $Y'_2\cap Z$ are
unions of irreducible components of $Y'\cap Z$ (for any
irreducible component $C$ of $Y\cap X'$ meets $Z$ properly in $X'$);
Moreover, their intersection has codimension $1$ in $Y'_1\cap Z$ and 
$Y'_2\cap Z$, by the first step of the proof. It follows that 
$Y'_1\cap Y'_2$ has codimension $1$ in both $Y'_1$ and $Y'_2$.
\end{proof}

\section{Singularities of orbit closures}

We begin by recalling the notion of rational singularities, see
e.g. \cite{KKMS} p.~50.

Let $Y$ be a variety. Choose a resolution of singularities
$\varphi:Z\to Y$, that is, $Z$ is smooth and $\varphi$ is proper and 
birational. Then the sheaves $R^i\varphi_*\cO_Z$ ($i\geq 0$) are
independent of the choice of $Z$. The singularities of $Y$ are
rational if $R^i\varphi_*\cO_Z=0$ for all $i\geq 1$ and
$\varphi_*\cO_Z=\cO_Y$; the latter condition is equivalent to
normality of $Y$. Varieties with rational singularities are
Cohen-Macaulay.

Let now $X$ be a spherical variety and $Y$ a $B$-stable subvariety. If
$Y$ is $G$-stable, then its singularities are rational, see e.g.
\cite{BI}. But this does not extend to arbitrary $Y$: generalizing
Example 1 in Section 1, we shall construct examples of $B$-orbit
closures of arbitrary dimension but of depth 1 at some points. In
particular, such orbit closures are neither normal nor Cohen-Macaulay.
\medskip

\noindent
{\sl Example 4.} Let $X$ be the space of unordered pairs $\{p,q\}$ of
distinct points in projective space $\mP^n$. The group 
$G={\rm GL}(n+1)$ acts transitively on $X$; one checks that $X$ is
spherical of rank 1. Let $\mP^m$ be a proper linear subspace of
$\mP^n$ of positive dimension $m$. Consider the space 
$$
m=\{\{p,q\}\in X~\vert~p\in\mP^m~{\rm or}~q\in\mP^m\},
$$
a subvariety of $X$ of codimension $n-m$. The stabilizer $P_m$ of
$\mP^m$ in $G$, a maximal parabolic subgroup, stabilizes $Y_m$ as
well; in fact, $Y_m$ contains an open $P_m$-orbit (the subset of all
$\{p,q\}$ such that $p\in\mP^m$ but $q\in\mP^n - \mP^m$) and its
complement
$$
Y'_m=\{\{p,q\}~\vert~p,q\in\mP^m,p\neq q\}
$$ 
is a unique $P_m$-orbit of codimension $n-m$ in $Y_m$. Thus, $Y_m$ is
the closure of a $B$-orbit; one checks that $r(Y_m)=0$ and
$r(Y'_m)=1$.

The map
$$
\begin{matrix}
\nu:&\mP^m\times\mP^n&\to& Y_m\cr
&(p,q)&\mapsto&\{p,q\}\cr
\end{matrix}
$$
is an isomorphism over the open $P_m$-orbit, but has degree $2$ over
$Y'_m$. Thus, $\nu$ is the normalization of $Y_m$, and the latter is not
normal. Moreover, $Y'_m$ is the singular locus of $Y_m$.

Observe that $Y_{n-1}$ is Cohen-Macaulay, as a divisor in $X$ (for
$n=2$ and $m=1$, we recover Example 1 in Section 1.) But if
$m<n-1$, then $Y_m$ has depth 1 along $Y'_m$ by Serre's criterion, see 
\cite{E} 18.3. In particular, $Y_m$ is not Cohen-Macaulay.

Let $\alpha_1,\ldots,\alpha_n$ be the simple roots of $G$. Then 
$P_{\alpha_m}Y_m=Y_{m+1}$, and $\alpha_m$ is the unique simple root
raising $Y_m$. The corresponding edge in $\Gamma(X)$ is simple, except
for $m=n-1$. Thus, $Y_m$ is the source of a unique oriented path with
target $X$, and the top edge of this path is double. In particular,
$Y_m$ is not multiplicity-free.
\medskip

Such examples of bad singularities do not occur for multiplicity-free
orbit closures:

\begin{theorem}\label{normal}
Let $Y$ be a multiplicity-free $B$-orbit closure in a spherical
$G$-variety $X$. If no simple normal subgroup of $G$ of type $G_2$,
$F_4$ or $E_8$ fixes points of $X$, then the singularities of $Y$ are
rational.
\end{theorem}

\begin{proof}
We begin with a reduction to the case where no simple normal
subgroup of $G$ fixes points of $X$. For this, we may assume that $G$
is the direct product of a torus with a family of simple, simply
connected subgroups; let $\Gamma$ be one of them. If $\Gamma$ is not
of type $G_2$, $F_4$ or $E_8$, then there exists a simple, simply
connected group $\tilde\Gamma$ together with a maximal proper
parabolic subgroup $\tilde P$ such that a Levi subgroup $\tilde L$ has
the same adjoint group as $\Gamma$ (indeed, add an edge to the Dynkin
diagram of $\Gamma$ to obtain that of $\tilde\Gamma$.) Then $\tilde L$
is the quotient of $\Gamma\times\mC^*$ by a finite central subgroup
$F$. We may assume moreover that $\mC^*$ maps injectively to 
$\tilde L$, that is, $\mC^*\cap F$ is trivial. Then the first
projection $p_1:F\to\Gamma$ is injective. 

We claim that the second projection $p_2:F\to\mC^*$ is injective as
well. Indeed, as $\tilde\Gamma$ is simply connected, its Picard group
is trivial; as some open subset of $\tilde\Gamma$ is the direct
product of $\tilde L$ with an affine space, the Picard group of
$\tilde L$ is trivial as well. But $\tilde L=(\Gamma\times\mC^*)/F$ is
the total space of the line bundle over $\Gamma/p_1(F)$ associated
with the character $p_2$ of $p_1(F)\cong F$, minus the zero
section. Thus, ${\rm Pic}(\tilde L)$ is the quotient of 
${\rm Pic}(\Gamma/p_1(F))$ by the class of that line bundle. Moreover,
${\rm Pic}(\Gamma/p_1(F))$ is isomorphic to the character group of
$F$, as $\Gamma$ is simply connected. Therefore $p_2$ generates the
character group of $F$. Since $F$ is abelian, the claim follows.

By that claim, $\Gamma\cap F$ is trivial; thus, $\Gamma$ embeds into
$\tilde L$ as its derived subgroup. We shall treat $p_2:F\to\mC^*$ as
an inclusion, which defines an action of $F$ on $\mC^*$. On the other
hand, $F$ acts on $X$ via $p_2:F\to\Gamma$, and this action commutes
with that of the remaining factors of $G$. Thus,
$X\times^F\mC^*$ is a variety with an action of the product
$\Gamma\times^F\mC^*\cong \tilde L$ with the remaining factors of
$G$. This variety is spherical and fibers equivariantly over
$\mC^*/F\cong \mC^*$, with fiber $X$. Thus, 
we may assume that the action of $\Gamma$ on $X$ extends to an action
of $\tilde L$. Now the parabolically induced variety 
$\tilde\Gamma\times^{\tilde P}X$ contains $Y$ as a multiplicity-free
subvariety (Lemma \ref{induced}) but contains no fixed point of
$\tilde\Gamma$. Iterating this argument removes the fixed points of all
simple normal subgroups of $G$.

We now reduce to the case where $X$ is projective. For this, we use
embedding theory of spherical homogeneous spaces, see \cite{K0}. We
may assume that $X$ contains a unique closed $G$-orbit $Z$ (for $X$ is
the union of $G$-stable open subsets, each of which contains a unique
closed $G$-orbit.) Together with Lemma \ref{parab}, the assumption
that no simple factor of $G$ fixes points of $X$ amounts to: $P(Z)$
contains no simple factor of $G$. Let $\cD_Z(X)$ be the set of all
colors $D$ that contain $Z$; then we can find an equivariant
projective completion $\overline{X}$ of $X$ such that 
$\cD_{Z'}(X)\subseteq \cD_Z(X)$ for any $G$-orbit closure $Z'$ in
$\overline{X}$. By Lemma \ref{parab}, it follows that 
$P(Z')\subseteq P(Z)$, and that no simple factor of $G$ fixes points
of $\overline{X}$.

We next reduce to an affine situation, in the following standard
way. Choose an ample $G$-linearized line bundle $\cL$ over
$X$. Replacing $\cL$ by a positive power, we may assume that $\cL$ is
very ample and that $X$ is projectively normal in the corresponding
projective embedding. Let $\hat X$ be the affine cone over $X$. This
is a spherical variety under the group $\hat G=G\times\mC^*$, and the
origin $0$ is the unique fixed point of any simple normal subgroup 
of $\hat G$, since $[\hat G,\hat G]=[G,G]$. Moreover, the affine cone
$\hat Y$ over $Y$ is stable under the Borel subgroup $B\times\mC^*$ 
of $\hat G$, and is multiplicity-free. Thus, we may assume that $X$ is 
affine with a fixed point $0$, and we have to show that $Y$ has
rational singularities outside $0$.

By \cite{BI}, the $G$-variety $GY$ is spherical, with rational
singularities, so that we may assume that $GY=X$. We argue then by
induction on the codimension of $Y$ in $X$.

Let $N_G(Y)$ be the set of all $g\in G$ such that $gY=Y$. This
is a proper standard parabolic subgroup of $G$, acting on $Y$ by
automorphisms. Let 
$$
\varphi:Z\to Y
$$ 
be a $N_G(Y)$-equivariant resolution of singularities. Denote by
$\mC[Y]$ (resp. $\mC[Z]$) the algebra of regular functions on $Y$
(resp. $Z$). Then $\mC[Z]$ is a finite $\mC[Y]$-module. Moreover, we
have an exact sequence of $\mC[Y]$-modules
$$
0\to\mC[Y]\to\mC[Z]\to C\to 0
$$
where the support of $C$ is the non-normal locus $N$ of $Y$, by
Zariski's main theorem. Note that $N_G(Y)$ acts on $C$ compatibly with
its $\mC[PY]$-module structure. We first show that $C$ is
supported at $0$, that is, $Y$ is normal outside $0$.

Let $\alpha$ be a simple root raising $Y$ and let $P=P_{\alpha}$. Let
$$
f=f_{Y,\alpha}:P\times^B Y \to P/B
$$ 
be the fiber bundle with fiber the $B$-variety $Y$; let
$$
\pi=\pi_{Y,\alpha}:P\times^B Y\to PY
$$
be the natural morphism. Then the map
$$
\pi^*:\mC[PY]\to\mC[P\times^B Y]
$$
is injective, and makes $\mC[P\times^B Y]$ a finite
$\mC[PY]$-module. Since $Y$ is multiplicity-free, $\pi$ is birational
and $PY$ is multiplicity-free as well. By the induction assumption,
$PY$ is normal outside $0$. Therefore, the cokernel of $\pi^*$ is
supported at $0$, by Zariski's main theorem again.

The $B$-equivariant resolution $\varphi:Z\to Y$ induces a
$P$-equivariant resolution
$$
\rho:P\times^B Z\to P\times^B Y. 
$$
Composing with $\pi$, we obtain a $P$-equivariant birational morphism
$$
{\tilde\pi}:P\times^B Z\to PY.
$$
As above, the map 
$$
{\tilde\pi}^*:\mC[PY]\to\mC[P\times^B Z]
$$
is injective and its cokernel is supported at $0$. We shall treat
$\pi^*$ and $\tilde\pi^*$ as inclusions.

We have
$$
\mC[P\times^B Y]=H^0(P\times^B Y,\cO_{P\times^B Y})
=H^0(P/B,f_*\cO_{P\times^B Y}).
$$
Moreover, $f_*\cO_{P\times^B Y}$ is the $P$-linearized sheaf on $P/B$
associated with the (rational, infinite-dimensional) $B$-module
$$
H^0(f^{-1}(B/B),\cO_{P\times^B Y})=\mC[Y].
$$
We shall use the notation 
$$
f_*\cO_{P\times^B Y}=\underline{\mC[Y]}.
$$ 
Then 
$$
\mC[PY]\subseteq H^0(P/B,\underline{\mC[Y]})
\subseteq H^0(P/B,\underline{\mC[Z]})=\mC[P\times^B Z]
$$ 
and these $\mC[PY]$-modules coincide outside $0$.

Consider the exact sequence of $P$-linearized sheaves on $P/B$:
$$
0\to\underline{\mC[Y]}\to \underline{\mC[Z]}\to
\underline{C}\to 0.
$$
Since the restriction map $\mC[PY]\to\mC[Y]$ is surjective, the
$B$-module $\mC[Y]$ is the quotient of a rational $P$-module. Since
$P/B$ is a projective line, it follows that
$H^1(P/B,\underline{\mC[Y]})=0$. Thus, we have an exact sequence of
$\mC[PY]$-modules 
$$
0\to H^0(P/B,\underline{\mC[Y]})\to 
H^0(P/B,\underline{\mC[Z]})
\to H^0(P/B,\underline{C})\to 0.
$$
It follows that $H^0(P/B,\underline{C})$ is supported at $0$.
Now normality of $Y$ outside $0$ is a consequence of the following 

\begin{lemma}\label{seshadri}
Let $C$ be a finite $\mC[Y]$-module with a compatible action of
$N_G(Y)$, such that the $\mC[PY]$-module $H^0(P/B,\underline{C})$ is
supported at $0$ for any minimal parabolic subgroup $P$ that raises
$Y$. Then $C$ is supported at $0$.
\end{lemma}

\begin{proof}
Otherwise, choose an irreducible component $Y'\neq\{0\}$ of
the support of $C$. Let $I(Y')$ be the ideal of $Y'$ in
$\mC[Y]$. Define a submodule $C'$ of $C$ by 
$$
C'=\{c\in C~\vert~ I(Y')c=0\}.
$$ 
Observe that the support of $C'$ is $Y'$ (indeed, the ideal $I(Y')$ is
a minimal prime of the support of $C$; thus, this ideal 
is an associated prime of $C$.) Note that $N_G(Y)$ stabilizes $Y'$
and acts on $C'$. Moreover, $H^0(P/B,\underline{C'})$ is a
$\mC[PY']$-module supported at $0$ (as a $\mC[PY]$-submodule of 
$H^0(P/B,\underline{C})$.)

We claim that $Y'$ is $G$-stable. Otherwise, let $\alpha$ be 
a simple root raising $Y'$; then $\alpha$ raises $Y$.
Define as above the maps 
$$
f':P\times^B Y'\to P/B~{\rm and }~\pi':P\times^B Y'\to PY'.
$$
The $\mC[Y']$-module $C'$ with a compatible $B$-action induces a
$P$-linearized sheaf $\cC'$ on $P\times^B Y'$, and we have
$f'_*\cC'=\underline{C'}$ as $P$-linearized sheaves on $P/B$. It
follows that the $\mC[PY']$-module
$H^0(P\times^B Y',\cC') = H^0(P/B,\underline{C'})$
is supported at $0$. On the other hand, we have
$H^0(P\times^B Y',\cC') = H^0(PY',\pi'_*\cC')$.
Moreover, the map $\pi':P\times^B Y'\to PY'$ is generically finite (as
$P$ raises $Y'$), and the support of $\cC'$ is $P\times^B Y'$ (as the
support of $C'$ is $Y'$). Thus, the support of $\pi'_*\cC'$ is $PY'$,
and the same holds for the support of 
$H^0(PY',\pi'_*\cC')=H^0(P/B,\underline{C'})$. This contradicts the
assumption that $Y'\neq\{0\}$. The claim is proved.

Let $L$ be the Levi subgroup of $P$ containing $T$, then
$P/B=[L,L]/B\cap [L,L]$. Since $Y'$ is $G$-stable, it is not fixed
pointwise by $[L,L]$ (here we use the assumption that no simple
normal subgroup of $G$ fixes points of $X - \{0\}$.) Since $Y'$ is
affine, $[L,L]$ acts non trivially on $\mC[Y']$. Thus, 
we can find an eigenvector $f$ of $B\cap [L,L]$ in 
$\mC[Y']=\mC[PY']$ of positive weight with respect to the coroot
$\check\alpha$. Then $f(0)=0$, so that $f$ acts nilpotently on
$H^0(P/B,\underline{C'})$. But $f$ does not act nilpotently on $C'$,
for the support of this module is $Y'$. Therefore we can choose a
finite-dimensional $B\cap [L,L]$-submodule $M$ of $C'$ such that 
$f^n M\neq 0$ for any large integer $n$. For such $n$, all weights of
$\check\alpha$ in $f^n M$ are positive. It follows that 
$H^0([L,L]/B\cap [L,L],\underline{f^n M})\neq 0$. But
$$
H^0([L,L]/B\cap [L,L],\underline{f^n M})\subseteq
H^0(P/B,\underline{f^n C'})=f^n H^0(P/B,\underline{C'}).
$$
Since $H^0(P/B,\underline{C'})$ is supported at $0$, we have 
$f^n H^0(P/B,\underline{C'})=0$ for large $n$, a contradiction. 
\end{proof}

Next we fix $i\geq 1$ and consider $R^i\varphi_*\cO_Z$, a
$N_G(Y)$-linearized coherent sheaf on $Y$. Since $Y$ is affine, this
sheaf is associated with the $\mC[Y]$-module $H^i(Z,\cO_Z)$ endowed
with a compatible action of $N_G(Y)$. We claim that the
$\mC[PY]$-module $H^0(P/B,\underline{H^i(Z,\cO_Z)})$ is supported at
$0$.

For this, note that the map $\tilde\pi:P\times^B Z\to PY$
is a resolution of singularities. By the induction assumption, $PY$
has rational singularities outside $0$; thus, the $\mC[PY]$-modules 
$H^q(P\times^B Z,\cO_{P\times^B Z})$ are supported at $0$, for all
$q\geq 1$. Moreover, $\tilde\pi = \pi\circ\rho$ (recall that
$\rho:P\times^B Z\to P\times^B Y$ denotes the $P$-equivariant
extension of $\varphi$.) And the fibers of $\pi:P\times^B Y\to PY$ 
identify to closed subsets of projective line, as the map 
$(\pi,f):P\times^B Y\to PY\times P/B$ 
is a closed immersion. Thus, $H^p(P\times^B Y,\cF)=0$ for any 
$p\geq 2$ and for any coherent sheaf $\cF$ on $P\times^B Y$. It
follows that the Leray spectral sequence 
$$
H^p(P\times^B Y,R^q\rho_*\cO_{P\times^B Z})\Rightarrow 
H^{p+q}(P\times^B Z,\cO_{P\times^B Z})
$$
degenerates at $E_2$: then 
$H^0(P\times^B Y,R^q\rho_*\cO_{P\times^B Z})$
is a quotient of $H^q(P\times^B Z,\cO_{P\times^B Z})$. In particular,
the $\mC[PY]$-module $H^0(P\times^B Y,R^i\rho_*\cO_{P\times^B Z})$ is
supported at $0$. Moreover, $R^i\rho_*\cO_{P\times^B Z}$ is the
$P$-linearized sheaf on $P\times^B Y$ associated with the
$B$-linearized sheaf $R^i\varphi_*\cO_Z$. Thus, 
$$
H^0(P\times^B Y,R^i\rho_*\cO_{P\times^B Z})=
H^0(P/B,\underline{H^i(Z,\cO_Z)}).
$$
This proves the claim. 

By Lemma \ref{seshadri}, it follows that the $\mC[Y]$-module
$H^i(Z,\cO_Z)$ is supported at $0$. Thus, $Y$ has rational
singularities outside $0$. 
\end{proof}

Combining Theorem \ref{normal} with Corollary \ref{simple}, we obtain
examples of spherical varieties where all $B$-orbit closures have
rational singularities, e.g., all embeddings of the symmetric spaces
listed at the end of Section 1. Here are other examples, of geometric
interest.
\medskip

\noindent
{\sl Example 5}. Let $\cF_n$ be the variety of all complete flags in
$\mC^n$. Consider the variety $X=\mP^{n-1}\times\cF_n$ endowed
with the diagonal action of $G={\rm GL}(n)$. Then $X$ is spherical,
see e.g. \cite{MWZ}. Clearly, the isotropy group of any point of $X$ 
is contained in a Borel subgroup of $G$; thus, by Lemma \ref{length},
all $B$-orbit closures in $X$ are multiplicity-free. Applying Theorem
\ref{normal}, it follows that their singularities are
rational. Therefore all ${\rm GL}(n)$-orbit closures in 
$\mP^{n-1}\times\cF_n\times\cF_n$ have rational singularities as well.
\medskip

\noindent
{\sl Example 6.} Let $p$, $q$, $n$ be positive integers such that 
$p\leq q\leq n$. Let $\cG_{n,p}$ be the Grassmanian variety of all 
$p$-dimensional linear subspaces of $\mC^n$. Consider the variety
$X=\cG_{n,p}\times\cG_{n,q}$ endowed with the diagonal action of 
$G={\rm GL}(n)$. By \cite{L}, $X$ is spherical (see also \cite{MWZ}.)

We claim that all edges of $\Gamma(X)$ are simple. Thus, the 
singularities of all $B$-orbit closures in $X$ are rational, 
and the same holds for closures of ${\rm  GL}(n)$-orbits in
$\cG_{n,p}\times\cG_{n,q}\times\cF_n$. 

To prove the claim, consider a point $(E,F)$ in the open 
$G$-orbit in $X$. Let $r=\dim(E\cap F)$, then 
$r={\rm max}(p+q-n,0)$. We can 
choose a basis $(v_1,\ldots,v_n)$ of $\mC^n$ such that $E\cap F$ 
(resp. $E$; $F$) is spanned by $v_1,\ldots,v_r$ 
(resp. $v_1,\ldots,v_p$; $v_1,\ldots,v_r,v_{p+1},\ldots,v_{p+q-r}$). 
Then, in the corresponding decomposition 
$$
\mC^n=\mC^r\oplus\mC^{p-r}\oplus\mC^{q-r}\oplus\mC^{n-p-q+r},
$$
the isotropy group of $(E,F)$ in $G$ consists of the following 
block matrices:
$$
\left(\begin{matrix}
*&*&*&*\cr
0&*&0&*\cr
0&0&*&*\cr
0&0&0&*\cr
\end{matrix}\right).
$$
Thus, the orbit $G/G_{(E,F)}$ is induced from 
${\rm GL}(n-r)/{\rm GL}(p-r)\times{\rm GL}(q-r)$. Now the claim
follows from Lemma \ref{induced} together with Corollary \ref{simple}.
\medskip

\noindent
{\sl Remark}. The varieties $\mP^{n-1}\times\cF_n\times\cF_n$ and
$\cG_{n,p}\times\cG_{n,q}\times\cF_n$
are examples of ``multiple flag varieties of finite type'' in the sense 
of \cite{MWZ}. There these varieties are classified for $G=GL(n)$. Do
all orbit closures in such varieties have rational singularities ? 
\medskip

\noindent
{\sl Example 7}. Let $M_{m,n}$ be the space of all 
$m\times n$ matrices. This is a spherical variety for the action of
$G={\rm GL}(m)\times{\rm GL}(n)$ by left and right
multiplication. Arguing as in Example 6, one checks that all $B$-orbit
closures in $M_{m,n}$ are multiplicity-free (in fact, any
$Y\in\cB(M_{m,n})$ satisfies $r(Y)=r(GY)$). Hence they have rational
singularities, by Theorem \ref{normal}.

The same result holds for the natural action of ${\rm GL}(n)$ on the
space of antisymmetric $n\times n$ matrices; but it fails in the case
of symmetric $n\times n$ matrices, if  $n\geq 3$. Indeed, the subset
$$
a_{11}=\left\vert
\begin{matrix}
a_{11}& a_{12}& a_{13}\cr
a_{12}& a_{22}& a_{23}\cr
a_{13}& a_{23}& a_{33}\cr
\end{matrix}
\right\vert=0
$$
is irreducible, stable under the standard Borel subgroup of $G$, and
singular along its divisor $(a_{11}=a_{12}=a_{13}=0)$.
\medskip

\begin{theorem}\label{reduced}
Let $X$ be a regular $G$-variety, let $Y$ be a
multiplicity-free $B$-orbit closure in $X$ such that $GY=X$, and let
$X'$ be a $G$-orbit closure in $X$, transversal intersection of the
boundary divisors $D_1,\ldots,D_r$. Then the singularities of $Y$ are
rational, and the scheme-theoretical intersection $Y\cap X'$ is
reduced. Moreover, for any $y\in Y\cap X'$, local equations
of $D_1,\ldots,D_r$ at $y$ are a  regular sequence in $\cO_{Y,y}$.
\end{theorem}

\begin{proof}
For rationality of singularities of $Y$, it is enough to check that
$X$ satisfies the assumption of Theorem \ref{normal}. We may assume
that $G$ acts effectively on $X$. If a simple normal subgroup $\Gamma$ 
of $G$ fixes points of $X$, let $X'$ be a component of the fixed
point set. Then $X'$ is $G$-stable: it is the closure of some orbit
$Gx$. Since $X$ is regular, the normal space $T_x(X)/T_x(Gx)$ 
is a direct sum of $\Gamma$-invariant lines. Since $\Gamma$ is simple 
and fixes pointwise $Gx$, it fixes pointwise $T_x(X)$ as well. 
It follows that $\Gamma$ fixes pointwise $X$, a contradiction.

For the remaining assertions, observe that the local equations of
$D_1,\ldots,D_r$ at any point $x\in X'$ are a regular sequence in
$\cO_{X,x}$. Moreover, as noted above, the scheme-theoretical
intersection $Y\cap X'$ is equidimensional of codimension $r$, and
generically reduced. Since $Y$ is Cohen-Macaulay, then $Y\cap X'$ is
reduced, and the local equations of $D_1,\ldots,D_r$ at any point 
$y\in Y\cap X'$ are a regular sequence in $\cO_{Y,y}$.
\end{proof}

We now apply these results to orbit closures in flag varieties.
For this, we recall a construction from \cite{B1} 1.5.
Let $G/H$ be a spherical homogeneous space, then $H$ acts on the flag
variety $G/B$ with only finitely many orbits. Let $V$ be a $H$-orbit 
closure in $G/B$ and let $\hat V$ be the corresponding $B$-orbit 
closure in $G/H$. Choose a complete regular embedding $X$ of
$G/H$ and let $Y$ be the closure of $\hat V$ in $X$. Then
$Y\in\cB(X)$ and $GY=X$. Consider the natural morphism 
$$
\pi:G\times^B Y\to X
$$ 
and the projection 
$$
f:G\times^B Y\to G/B.
$$
The fibers of $\pi$ identify to closed subschemes of $G/B$ via
$f_*$. Let $x$ be the image in $X$ of the base point of $G/H$, then
$\pi^{-1}(x)$ identifies to $V$. On the other hand, let $Z$ be a
closed $G$-orbit in $X$ with $B$-fixed point $z$, then the set
$f(\pi^{-1}(z))$ equals
$$
V_0=\bigcup_{w\in W(Y)}\overline{Bw_0wB}/B
$$
where $w_0$ denotes the longest element of $W$. Moreover, we have in
the integral cohomology ring of $G/B$:
$$
[V]=\sum_{w\in W(Y)} d(Y,w)[\overline{Bw_0wB}/B].
$$

Now Theorem \ref{connected} and Proposition \ref{constant} imply the
following

\begin{corollary}\label{flag}
Notation being as above, $V_0$ is connected in codimension $1$. If
moreover $G$ is simply-laced, then $[V]=2^{\ell_N(\gamma)}[V_0]$ where
$\gamma$ is any oriented path in $\Gamma(X)$ joining $Y$ to $X$.
\end{corollary}

We shall call $V$ multiplicity-free if $Y$ is. Equivalently, the 
cohomology class of $V$ decomposes as a sum of Schubert classes with
coefficients $0$ or $1$.

Note that any multiplicity-free $H$-orbit closure $V$ is irreducible,
even if $H$ is not connected. Indeed, $H$ acts transitively on the set
of all irreducible components of $V$, so that any two such components
have the same cohomology class; but the class of $V$ is indivisible in
the integral cohomology of $G/B$.

\begin{theorem}\label{flat}
Let $G/H$ be a spherical homogeneous space, and $V$ a 
multiplicity-free $H$-orbit closure in $G/B$. Then the singularities
of $V$ are rational.

Moreover, let $X$ be a complete regular embedding of $G/H$ 
and let $Y$ be the $B$-orbit closure in $X$ associated with $V$, 
then the natural morphism $\pi:G\times^B Y\to X$ is flat, and its 
fibers are reduced. 

As a consequence, the fibers of $\pi$ realize a degeneration of $V$ to 
the reduced subscheme $V_0$ of $G/B$. 
\end{theorem}

\begin{proof} Note that the singularities of $Y$ are rational by
Theorem \ref{reduced}; thus, the same holds for $\hat V=Y\cap G/H$. 
Let $\varphi:Z\to \hat V$ be a resolution of singularities; consider the
quotient map $q_H:G\to G/H$, the preimage $V'=q_H^{-1}(\hat V)$ in
$G$, and the fiber product $Z'=Z\times_{\hat V} V'$. Then $V'$ is
smooth, since $Z$ and $q_H$ are; the projection $\varphi':Z'\to V'$ is
proper and birational, since $\varphi$ is; and
$R^i\varphi_*\cO_{Z'}=0$ for $i\geq 1$, since cohomology commutes with
flat base extension. Therefore the singularities of $V'$ are rational.
Now $V'=q_B^{-1}(V)$ and $q_B$ is a locally trivial fibration, so that
the singularities of $V$ are rational as well. 

For the second assertion, we identify $Y$ to its image $B\times^B Y$ 
in $G\times^B Y$. Since $\pi$ is $G$-equivariant, it is enough to check 
the statement at $y\in Y$. Let $D_1,\ldots,D_r$ be the boundary divisors 
containing $y$, with local equations $f_1,\ldots,f_r$ in $\cO_{X,y}$. 
It follows from Theorem \ref{reduced} that the pull-backs 
$\pi^*f_1,\ldots,\pi^*f_r$ are a regular sequence in 
$\cO_{G\times^B Y,y}$ and generate the ideal of $\pi^{-1}(Gy)$. 
Moreover, the restriction of $\pi$ to $\pi^{-1}(Gy)$ is flat with
reduced fibers, as $\pi$ is $G$-equivariant. Now we conclude
by a local flatness criterion, see \cite{E} Corollary 6.9.
\end{proof}

A direct consequence is the following

\begin{corollary}\label{surj}
Consider a spherical homogeneous space $G/H$, a multiplicity-free
$H$-orbit closure $V$ in $G/B$ and an effective line bundle $L$ on
$G/B$. Then the restriction map $H^0(G/B,L)\to H^0(V,L)$ is
surjective, and $H^i(V,L)=0$ for all $i\geq 1$.
\end{corollary}

Indeed, this holds with $V$ replaced by $V_0$, a union of Schubert
varieties (see \cite{MR}.) The result follows by semicontinuity of
cohomology in a flat family.

We now obtain a partial converse to Corollary \ref{surj}:

\begin{proposition}\label{nonsurj}
Let $G/H$ be a spherical homogeneous space, let $V$ be a $H$-orbit
closure in $G/B$ and let $Y$ be the corresponding $B$-orbit closure in
$G/H$. If $Y$ is the source of a double edge of $\Gamma(G/H)$, then 
there exists an effective line bundle $L$ on $G/B$ such that the
restriction $H^0(G/B,L)\to H^0(V,L)$ is not surjective.
\end{proposition}

\begin{proof}
Let $\alpha$ be the label of a double edge with source $Y$. Denote by
$p:G/B\to G/P_{\alpha}$ the natural map and by $p_V:V\to \pi(V)$
its restriction to $V$; then $p$ is a projective line bundle, and
$p_V$ is generically finite of degree $2$. Choose an ample line
bundle $L$ on $G/P_{\alpha}$; then $p^*L$ is an effective line
bundle on $G/B$. Now our assertion is a direct consequence of the
following claim: the restriction map
$$
r_n:H^0(p^{-1}p(V),p^*(L^{\otimes n}))\to 
H^0(V,p^*(L^{\otimes n}))
$$
is not surjective for large $n$. To check this, note that
$H^0(p^{-1}p(V),p^*(L^{\otimes n}))=H^0(p(V),L^{\otimes n})$ 
and that
$H^0(V,p^*(L^{\otimes n}))=
H^0(p(V),L^{\otimes n}\otimes p_{V*}\cO_V)$, by the projection
formula. Thus, $r_n$ identifies with the map
$$
H^0(p(V),L^{\otimes n})\to 
H^0(p(V),L^{\otimes n}\otimes p_{V*}\cO_V)
$$
defined by the inclusion of $\cO_{p(V)}$ into $p_{V*}\cO_V$. Since
$p_V$ has degree $2$, the quotient $\cF=p_{V*}\cO_V/\cO_{p(V)}$
has rank $1$ as a sheaf of $\cO_{p(V)}$-modules. Moreover, since $L$
is ample, the cokernel of $r_n$ is isomorphic to 
$H^0(p(V),\cF\otimes L^{\otimes n})$ for large $n$. This proves the
claim. 
\end{proof}

\section{Orbit closures of maximal rank}

Let $\cB(X)_{max}$ be the set of all $Y\in\cB(X)$ such that 
$r(Y)=r(X)$, that is, the set of all $B$-orbit closures of maximal
rank. Recall that all such orbit closures are multiplicity-free and
meet the open $G$-orbit. Here is another characterization of them.

\begin{proposition}\label{max}
(i) For any $Y\in\cB(X)_{max}$ and $w\in W(Y)$, we have:
$BwY^0=X^0$ and $w^{-1}\in W^{\Delta(X)}$. Moreover, $W(Y)$ is
disjoint from all $W(Y')$ where $Y'\in\cB(X)$ and $Y'\neq Y$.

\noindent
(ii) Conversely, if $Y\in\cB(X)$ and there exists $w\in W$ such that
$BwY^0=X^0$, then $Y$ has maximal rank. If moreover 
$w^{-1}\in W^{\Delta(X)}$, then $w\in W(Y)$, and $\Delta(Y)$
consists of those $\alpha\in\Delta$ such that $w(\alpha)\in\Delta(X)$.
\end{proposition}

\begin{proof}
(i) We prove that $BwY^0=X^0$ by induction over $\ell(w)$, the
case where $\ell(w)=0$ being evident. If $\ell(w)\geq 1$, we can write 
$w=w's_{\alpha}$ for some simple root $\alpha$ and some $w'\in W$ such
that $\ell(w')=\ell(w)-1$; then $BwB=Bw'Bs_{\alpha}B$. Then
$X=\overline{BwY}=\overline{Bw'P_{\alpha}Y}$. Since
$\ell(w)={\rm codim}_X(Y)$, it follows that $\alpha$ raises $Y$ and
that $w'\in W(P_{\alpha}Y)$. Because $Y$ has maximal rank,
$P_{\alpha}Y^0$ consists of two $B$-orbits, both of maximal rank. But
$P_{\alpha}Y^0=Y^0\cup Bs_{\alpha}Y^0$ so that $Bs_{\alpha}Y^0$ is a
unique $B$-orbit of maximal rank and of codimension $\ell(w')$ in
$X$. By the induction assumption, we have $Bw'Bs_{\alpha}Y^0=X^0$,
that is, $BwY^0=X^0$. If moreover $w\in W(Y')$ for some $Y'\in\cB(X)$,
then a similar induction shows that $Y'=Y$.

If $w^{-1}\notin W^{\Delta(X)}$ then there exists $\beta\in\Delta(X)$
such that $\ell(s_{\beta}w)=\ell(w)-1$. Thus,
$BwB=Bs_{\beta}Bs_{\beta}wB$, so that $s_{\beta}Bs_{\beta}wY^0$ is
contained in $X^0$. But $s_{\beta}X^0=X^0$; therefore, 
$Bs_{\beta}wY^0 = X^0$, and $\overline{Bs_{\beta}wY}=X$. It follows
that ${\rm codim}_X(Y)\leq\ell(s_{\beta}w)=\ell(w)-1$, a
contradiction.

(ii) Let ${\dot w}$ be a representative of $w$ in the normalizer of
$T$. By assumption, the map 
$$
\begin{matrix}U\times Y^0&\to&X^0\cr
(u,y)&\mapsto&u{\dot w}y\cr
\end{matrix}
$$
is surjective. Thus, it induces an injective homomorphism from the
ring $\mC[X^0]$ of regular functions on $X^0$, to $\mC[U\times Y^0]$. 
The group of invertible regular functions $\mC[X^0]^*$ is mapped into 
$\mC[U\times Y^0]^*=\mC[Y^0]^*$. Quotienting by
$\mC^*$ and taking ranks, we obtain $r(X)\leq r(Y)$ by Lemma
\ref{prelim}, whence $r(Y)=r(X)$. 

If moreover $w^{-1}\in W^{\Delta(X)}$, we show that $w\in W(Y)$ by
induction over $\ell(w)$; we may assume that $w\neq 1$. Then we can
write $w=w's_{\alpha}$ where $w'\in W$, $\alpha\in\Delta$ and
$\ell(w)=\ell(w')+1$. It follows that $w(\alpha)\in\Phi^-$. 

We begin by checking that $s_{\alpha} Y^0\neq Y^0$. Otherwise, by
Lemma \ref{prelim}, there exists $y\in Y^0$ fixed by
$[L_{\alpha},L_{\alpha}]$. Thus, ${\dot w}y\in X^0$ is fixed by
$w[L_{\alpha},L_{\alpha}]w^{-1}$. Since the unipotent radical of $P(X)$
acts freely on $X^0$ by Lemma \ref{parab}, it follows that
$w(\alpha)\in\Phi_{\Delta(X)}$. Then 
$\alpha\in\Delta\cap w^{-1}(\Phi_{\Delta(X)}^-)$
which contradicts the assumption that $w^{-1}\in W^{\Delta(X)}$. 

As above, it follows that $Bs_{\alpha}Y^0$ is a $B$-orbit of maximal
rank and of dimension $\dim(Y)+1$; moreover,
$Bw'Bs_{\alpha}Y^0=X^0$. We can write $w'=uv$ where 
$u\in W_{\Delta(X)}$, $v^{-1}\in W^{\Delta(X)}$, and
$\ell(w')=\ell(u)+\ell(v)$. Thus,
$BwB=BuBvBs_{\alpha}B$, and $BvBs_{\alpha}Y^0=X^0$ as
$u^{-1}X^0=X^0$. By the induction assumption, 
$v\in W(\overline{Bs_{\alpha}Y})$. Moreover,
$\ell(vs_{\alpha})=\ell(v)+1$, for $w=uvs_{\alpha}$ and
$\ell(w)=\ell(u)+\ell(v)+1$. It follows that $vs_{\alpha}\in W(Y)$;
in particular, $s_{\alpha}v^{-1}\in W^{\Delta(X)}$. But
$w^{-1}=s_{\alpha}v^{-1}u^{-1}$ is in $W^{\Delta(X)}$ as well. Thus,
$u=1$ and $w^{-1}\in W(Y)$.

Let $\alpha$ be a simple root of $Y$. Then we see as above that
$w(\alpha)\in\Phi_{\Delta(X)}$. We have 
$ws_{\alpha}=s_{w(\alpha)}w$ with $s_{w(\alpha)}\in W_{\Delta(X)}$ and
$w^{-1}\in W^{\Delta(X)}$. Thus,
$\ell(ws_{\alpha})=\ell(s_{w(\alpha)})+\ell(w)$ 
which forces $w(\alpha)\in\Phi^+$ (as $\ell(s_{\alpha}w)=\ell(w)+1$)
and $w(\alpha)\in\Delta$ (as $\ell(s_{w(\alpha)})=1$.) We conclude
that $w(\alpha)$ is a simple root of $X$. 

Conversely, let $\alpha\in\Delta$ such that $w(\alpha)$ is a simple
root of $X$. Then $\ell(ws_{\alpha})=\ell(w)+1$, whence 
$$
BwBs_{\alpha}Y^0=Bws_{\alpha}Y^0=Bs_{w(\alpha)}wY^0
=Bs_{w(\alpha)}BwY^0=Bs_{w(\alpha)}X^0=X^0.
$$
Let $\cO$ be a $B$-orbit in $Bs_{\alpha}Y^0$. Then $Bw\cO=X^0$. By
(i), we have $\cO=Y^0$, whence $s_{\alpha}Y^0=Y^0$ and
$\alpha\in\Delta(Y)$.
\end{proof}

This preliminary result, combined with those of Section 2, implies a
structure theorem for orbits of maximal rank and their closures in
regular varieties:

\begin{theorem}\label{local}
Let $X$ be a complete regular $G$-variety, $Y\in\cB(X)_{max}$ and
$w\in W(Y)$. Choose a ``slice'' $S_{Y,w}$ as in Proposition
\ref{slice}, so that the product map 
$$
(U\cap w^{-1}R_u(P)w)\times w^{-1}S_{Y,w}\to Y\cap w^{-1}X_0
$$ 
is an isomorphism. Then $w^{-1}S_{Y,w}$ is fixed pointwise by
$[L(Y),L(Y)]$. Moreover, $Y\cap w^{-1}X_0$ is $P(Y)$-stable and meets
each $G$-orbit along a unique $B$-orbit, of maximal rank in this
$G$-orbit. In particular, there exists $y\in Y^0$ fixed by
$[L(Y),L(Y)]$ such that the product map 
$(U\cap w^{-1}R_u(P)w)\times Ty\to Y^0$ 
is an isomorphism.

As a consequence, for each $G$-orbit closure $X'$ in $X$, all
irreducible components of $Y\cap X'$ have maximal rank
in $X'$. Moreover, a given $Y'\in\cB(X')$ is an irreducible component
of $Y\cap X'$ if and only if $W(Y')$ is contained in $W(Y)$.
\end{theorem}

\begin{proof}
With notation as in Section 2, recall that 
$$
w^{-1}S_{Y,w}=Y\cap (U^-\cap w^{-1}Uw)w^{-1}S
$$
where $S$ is fixed pointwise by $[L(X),L(X)]$. Now Proposition
\ref{max} implies that $[L(Y),L(Y)]$ fixes pointwise $S$ and
normalizes $U^-\cap w^{-1}Uw$. Thus, $[L(Y),L(Y)]$ stabilizes
$w^{-1}S_{Y,w}$. Moreover, intersecting that space with those boundary
divisors that contain a given closed $G$-orbit, we obtain
$[L(Y),L(Y)]$-stable hypersurfaces meeting transversally at a fixed
point. Arguing as in the proof of Theorem \ref{reduced}, it follows
that $[L(Y),L(Y)]$ fixes pointwise $w^{-1}S_{Y,w}$.

By Proposition \ref{slice}, $w^{-1}S_{Y,w}$ meets each
$G$-orbit along a unique $T$-orbit. As a consequence, the intersection
of $Y\cap w^{-1}X_0$ with each $G$-orbit is contained in a unique
$B$-orbit. We apply this to $GY^0$, the open $G$-orbit in $X$. Since 
$Y\cap w^{-1}X_0\cap GY^0=Y\cap w^{-1}X^0$ equals $Y^0$ by Proposition
\ref{max}, we see that the product map
$$
(U\cap w^{-1}R_u(P)w)\times (w^{-1}S_{Y,w}\cap Y^0)\to Y^0
$$
is an isomorphism. Moreover, $w^{-1}S_{Y,w}\cap Y^0$ is a unique
$T$-orbit of dimension equal to the rank of $X$.

It follows that each $U$-orbit in $Y^0$ is a unique orbit of 
$U\cap w^{-1}R_u(P)w$. Indeed, any $U$-orbit is isomorphic to some
affine space, and its projection to $w^{-1}S_{Y,w}\cap Y^0$ is a
morphism to a torus, hence is constant.

Choose $y_0\in Y^0$ and let $y\in Y\cap w^{-1}X_0$. Since 
$By_0=Y^0$ is dense in $Y\cap w^{-1}X_0$, we have
$\dim(Uy)\leq \dim(Uy_0)$. The latter equals
$\dim(U\cap w^{-1}R_u(P)w)$ by the previous step. Because
$U\cap w^{-1}R_u(P)w$ acts freely on $Y\cap w^{-1}X_0$, it follows
that $(U\cap w^{-1}R_u(P)w)y$ is open in $Uy$. But both are affine
spaces, so that they are equal. Thus, $Y\cap w^{-1}X_0$ is $B$-stable.
It is even $P(Y)$-stable, because $P(Y)\subseteq w^{-1}Pw$ by
Proposition \ref{max}. 

Since $w^{-1}S_{Y,w}$ meets each $G$-orbit along a unique $T$-orbit,
$Y\cap w^{-1}X_0$ meets each $G$-orbit along a unique $B$-orbit. Let
$y\in Y\cap w^{-1}X_0$, then $wBy\subseteq X_0$ and, therefore,
$wBy\subseteq(Gy)^0$. By Proposition \ref{max} again, we have
$r(By)=r(Gy)$.

The remaining assertions follow from Theorem \ref{indices} together
with Proposition \ref{max}.
\end{proof}

We now describe the intersections of $B$-orbit closures of maximal rank
with $G$-orbit closures, in terms of Knop's action of the Weyl group
$W$ on the set $\cB(X)$. This action can be defined as follows. 

Let $\alpha\in\Delta$ and $Y\in\cB(X)$, then $s_{\alpha}$ fixes $Y$,
except in the following cases:

\noindent
$\bullet$ Type $U$: $P_{\alpha}Y^0=Y^0\cup Z^0$ for $Z\in\cB(X)$ with
$r(Z)=r(Y)$. Then $s_{\alpha}$ exchanges $Y$ and $Z$.

\noindent
$\bullet$ Type $T$: $P_{\alpha}Y^0=Y^0\cup Y_-^0\cup Z^0$ for
$Z\in\cB(X)$ with $r(Y)=r(Y_-)=r(Z)-1$. Then $s_{\alpha}$ exchanges
$Y$ and $Y_-$.

By \cite[\S 4]{K2}, this defines indeed a $W$-action (that is, the
braid relations hold); moreover, $\cX(w(Y))=w(\cX(Y))$ for all $w\in W$.
In particular, this action preserves the rank. 

For $Y\in\cB(X)_{max}$ and $w\in W(Y)$, we have $w(Y)=X$. Thus, 
$\cB(X)_{max}$ is the $W$-orbit of $X$ in $\cB(X)$.

Let $W_{(X)}$ be the isotropy group of $X$; then $W_{(X)}$ acts on
$\cX(X)$. Observe that $W_{(X)}$ contains $W_{\Delta(X)}$. The latter
acts trivially on $\cX(X)$ by Lemma \ref{prelim}. In fact, $W_{(X)}$
stabilizes $\Phi_{\Delta(X)}$ (indeed, $\Phi_{\Delta(X)}$ consists of
all roots that are orthogonal to $\cX(X)$, if $X$ is non-degenerate in
the sense of \cite{K1}; and the general case reduces to that one, by
\cite{K1} \S 5.) 

The normalizer of $\Phi_{\Delta(X)}$ in $W$ is the semi-direct product
of $W_{\Delta(X)}$ with the normalizer of $\Delta(X)$. Therefore,
$W_{(X)}$ is the semi-direct product of $W_{\Delta(X)}$ with
$$
W_X=\{w\in W~\vert~w(X)=X~{\rm and}~w(\Delta(X))=\Delta(X)\}.
$$
The latter identifies to the image of $W_{(X)}$ in 
${\rm Aut}\,\cX(X)$, that is, to the ``Weyl group of $X$'', see
\cite{K2} Theorem 6.2. 

In fact, $W_X$ is the set of all $w\in W_{(X)}$ such that
$w(\rho)-\rho\in\cX(X)$, where $\rho$ denotes the half sum of positive
roots (see \cite{K} 6.5); we shall not need this result. 

Let
$$
W^{(X)}=\{w\in W~\vert~\ell(wu)\geq \ell(w)~\forall~u\in W_{(X)}\},
$$
the set of all elements of minimal length in their right
$W_{(X)}$-coset.

\begin{proposition}\label{minimal}
Notation being as above, we have
$$
W^{(X)}=\{w\in W^{\Delta(X)}~\vert~
\ell(wu)\geq \ell(w)~\forall~u\in W_X\},
$$
and, for any $w\in W$,
$$
W(w(X))=\{v\in W~\vert~ v^{-1}\in W^{(X)}\cap wW_{(X)}\}.
$$

As a consequence, all elements of minimal length in a given left
$W_{(X)}$-coset have the same length and are contained in a left
$W_X$-coset. Moreover, the subsets $W(Y)$, $Y\in\cB(X)_{max}$, are
exactly the subsets of all elements of minimal length in
a given left $W_{(X)}$-coset.

If moreover $X$ is regular, then we have for any $G$-orbit closure
$X'$ in $X$:
$$
w(X)\cap X'=\bigcup_{w'\in W^{(X)}\cap wW_{(X)}} w'(X').
$$
\end{proposition}

\begin{proof}
Clearly, $W^{(X)}$ is contained in $W^{\Delta(X)}$. And since $W_X$
stabilizes $\Delta(X)$, the set $W^{\Delta(X)}$ is stable under right
multiplication by $W_X$. This implies the first assertion.

Let $Y=w(X)$ and observe that ${\rm codim}_X(Y)\leq\ell(w)$ with
equality if and only if $w^{-1}\in W(Y)$ (indeed, a reduced
decomposition of $w$ defines a non-oriented path in $\Gamma(X)$ with
endpoints $Y$ and $X$).

Let $v\in W(Y)$. Since $v(Y)=X$, we have $v^{-1}\in wW_{(X)}$. 
Moreover, $\ell(v^{-1})=\ell(v)={\rm codim}_X(Y)\leq \ell(w)$. Since 
we can change $w$ in its right $W_{(X)}$-coset, it follows that 
$v^{-1}\in W^{(X)}$. 

Conversely, let $u\in W$ such that 
$u^{-1}\in W^{(X)}\cap wW_{(X)}$. Then $u(Y)=X$, whence 
$\ell(u)\geq\ell(v)$ and $u\in W_{(X)}v$. Since $u^{-1}\in W^{(X)}$,
this forces $\ell(u)=\ell(v)$ and then $u\in W(Y)$. This proves the
first assertion. Together with Theorem \ref{local}, this
implies the second assertion.
\end{proof}
\medskip

\noindent
{\sl Example 8}. Let ${\bf G}$ be a connected reductive
group. Consider the group $G={\bf G}\times{\bf G}$ acting on 
$X={\bf G}$ by $(x,y)\cdot z=xzy^{-1}$. Then $X$ is a spherical
homogeneous space: consider the Borel subgroup 
$B={\bf B}\times{\bf B}^-$ of $G$, where ${\bf B}$ and ${\bf B}^-$ are
opposed Borel subgroups of ${\bf G}$. With evident notation, the
$B$-orbits in $X$ are the ${\bf B}w{\bf B}^-$, $w\in{\bf W}$. 
This identifies $\cB(X)$ to ${\bf W}$. Moreover, all $B$-orbits have
maximal rank, and the Weyl group $W={\bf W}\times{\bf W}$ acts on 
${\bf W}$ by $(u,v)w=uwv^{-1}$. Thus, $\Delta(X)$ is empty, $W_{(X)}$
is the diagonal in ${\bf W}\times{\bf W}$, and ${\bf W}\times\{1\}$ is
a system of representatives of $W/W_{(X)}$. One checks that 
$$
W^{(X)}=\{(u,v)\in{\bf W}\times{\bf W}~\vert~
\ell(u)+\ell(v)=\ell(uv^{-1}\}.
$$
In particular, $(w,1)\in W^{(X)}$ for all $w\in{\bf W}$. Moreover,
$$
W^{(X)}\cap(w,1)W_{(X)}=\{(u,v)\in{\bf W}\times{\bf W}~\vert~
uv^{-1}=w~{\rm and}~\ell(u)+\ell(v)=\ell(w)\}.
$$
This identifies $W^{(X)}\cap(w,1)W_{(X)}$ to the set of all 
$u\in{\bf W}$ such that $u\preceq w$ for the right order on ${\bf W}$.
\medskip

\noindent
{\sl Remark}. Let $X$ be a complete regular $G$-variety, $Y$ a
$B$-orbit closure of maximal rank, and $X'$ a $G$-orbit closure in
$X$. Then the number of irreducible components of $Y\cap X'$ is at
most the order of $W_X$ by Proposition \ref{minimal}. If moreover $X$
has rank $1$, then $W_X$ is trivial or has order $2$, so that 
$Y\cap X'$ has at most $2$ components.

\medskip
 
Returning to an arbitrary spherical variety $X$, we shall deduce from
Proposition \ref{equal} the following

\begin{theorem}
The group $W_{(X)}$ is generated by reflections $s_{\alpha}$ where
$\alpha$ is a root such that $\alpha\in\Phi_{\Delta(X)}$ or that
$2\alpha\in\cX(X)$, and by products $s_{\alpha}s_{\beta}$ where
$\alpha$, $\beta$ are orthogonal roots such that 
$\alpha+\beta\in\cX(X)$.
\end{theorem}

\begin{proof}
Let $w\in W_{(X)}$. We choose a reduced decomposition
$w=s_{\alpha_{\ell}}\cdots s_{\alpha_2}s_{\alpha_1}$
and we argue by induction on $\ell$.

If $\alpha_1\in \Delta(X)$ then $s_{\alpha_1}$ is a reflection in
$W_{(X)}$, so that $s_{\alpha_{\ell}}\cdots s_{\alpha_2}\in
W_{(X)}$. Now we conclude by the induction assumption.

If $\alpha_1\notin\Delta(X)$ then $s_{\alpha_1}(X)$ has codimension 1
in $X$. Let $i$ be the largest integer such that 
${\rm codim}_Xs_{\alpha_i}\cdots s_{\alpha_1}(X)=i$. Let
$Y=s_{\alpha_i}\cdots s_{\alpha_1}(X)=i$, then $Y\in\cB(X)_{max}$ and 
$s_{\alpha_1}\cdots s_{\alpha_1}\in W(Y)$. 

If $P_{\alpha_{i+1}}Y=Y$ then $s_{\alpha_{i+1}}(Y)=Y$ by definition of
the $W$-action and maximality of $i$. Let 
$\alpha=s_{\alpha_1}\cdots s_{\alpha_i}(\alpha_{i+1})$.
Then $s_{\alpha}$ is a reflection of $W_{(X)}$, and
$w=s_{\alpha_{\ell}}\cdots s_{\alpha_{i+2}} s_{\alpha_i}\cdots
s_{\alpha_1}s_{\alpha}$. If $\alpha_{i+1}\in\Delta(Y)$, then
$\alpha\in\Delta(X)$ by Proposition \ref{max}. Otherwise,
$P_{\alpha_{i+1}}Y^0/R(P_{\alpha_{i+1}})$ is isomorphic to 
${\rm PGL}(2)/T$ or to ${\rm PGL}(2)/N$; it follows that
$2\alpha_{i+1}\in\cX(Y)$, and that $2\alpha\in\cX(X)$. Now we conclude
by the induction assumption.

If $P_{\alpha_{i+1}}Y\neq Y$ then $\alpha_{i+1}$ raises $Y$ to (say)
$Y'$.  Choose $u\in W(Y')$, then $\ell(u)=i-1$ and
$us_{\alpha_{i+1}}\in W(Y)$. Moreover, 
$us_{\alpha_{i+1}}s_{\alpha_i}\cdots s_{\alpha_1}\in W_{(X)}$. We have 
$w=vus_{\alpha_{i+1}}s_{\alpha_i}\cdots s_{\alpha_1}$ for some 
$v\in W_{(X)}$ such that $\ell(vu)=\ell-i-1$. Thus,
$\ell(v)\leq\ell(vu)+\ell(u)=\ell-2$. Therefore, we may assume that
there exist $Y\in\cB(X)_{max}$ and $w_1,w_2\in W(Y)$ such that
$w=w_2w_1^{-1}$. By Proposition \ref{neigh}, we may assume moreover
that $w_1$ and $w_2$ are neighbors. Then we conclude by Proposition
\ref{equal}. 

\end{proof}

As a direct consequence, we recover the following result of Knop, see
\cite{K1} and \cite{K2}.

\begin{corollary}
The image of $W_X$ in ${\rm Aut}\,\cX(X)$ is generated by
reflections.
\end{corollary}

\end{document}